\documentclass[12pt]{article}
\usepackage{amsmath}
\usepackage{amsfonts}
\usepackage{amssymb}
\usepackage{longtable}

\newtheorem{theorem}{Theorem}
\newtheorem{prop}{Proposition}
\newtheorem{lem}{Lemma}
\newtheorem{ex}{Example}
\newtheorem{tab}{Table}

\def\spa{\mathop\text{{\rm span}}\nolimits}

\def\pr{\mathop\text{\rm pr}\nolimits}
\def\tr{\mathop\text{\rm tr}\nolimits}
\def\str{\mathop\text{\rm str}\nolimits}
\def\Ric{\mathop\text{\rm Ric}\nolimits}

\def\ad{\mathop\text{\rm ad}\nolimits}

\def\GL{\text{\rm GL}}
\def\SL{\text{\rm SL}}
\def\SO{\text{\rm SO}}
\def\Ort{\text{\rm O}}
\def\Un{\text{\rm U}}
\def\Sp{\text{\rm Sp}}
\def\SU{\text{\rm SU}}

\def\Real{\mathbb{R}}
\def\Co{\mathbb{C}}

\def\g{\mathfrak{g}}
\def\h{\mathfrak{h}}

\def\z{\mathfrak{z}}

\def\so{\mathfrak{so}}
\def\osp{\mathfrak{osp}}
\def\hosp{\mathfrak{hosp}}
\def\sl{\mathfrak{sl}}

\def\sp{\mathfrak{sp}}
\def\gl{\mathfrak{gl}}
\def\su{\mathfrak{su}}
\def\u{\mathfrak{u}}
\def\hol{\mathfrak{hol}}
\def\q{\mathfrak{q}}
\def\pe{\mathfrak{pe}}

\def\spin{\mathfrak{spin}}
\def\z{\mathfrak{z}}

\def\Ga{\Gamma}
\def\na{\nabla}

\def\t{\tilde}
\def\ti{\thicksim}
\def\ga{\gamma}

\def\T{{\cal T}}
\def\M{{\cal M}}
\def\O{{\cal O}}
\def\Hols{{\mathcal Hol}}
\def\Hol{\mathop\text{\rm Hol}\nolimits}

\def\R{{\cal R}}
\def\RR{\bar{\cal R}}
\def\P{{\cal P}}

\begin{document}

\title{Irreducible holonomy algebras of Riemannian supermanifolds
}


\author{Anton S. Galaev}



\maketitle

\begin{abstract} Possible
irreducible holonomy algebras $\g\subset\osp(p,q|2m)$ of
Riemannian supermanifolds under the assumption that $\g$ is a
direct sum of simple Lie superalgebras of classical type and
possibly of a one-dimensional center are classified. This
generalizes the classical result of Marcel Berger about the
classification of irreducible holonomy algebras of
pseudo-Riemannian manifolds.

\end{abstract}

\tableofcontents

\section{Introduction} Possible irreducible holonomy algebras
(or, equivalently, connected irreducible holonomy groups) of
non-locally symmetric pseudo-Riemannian manifolds are classified
by M.~Berger in \cite{Ber}. These algebras are the following:
$\so(p,q)$, $\so(p,\Co)$,  $\u(r,s)$, $\su(r,s)$,  $\sp(r,s)$,
$\sp(r,s)\oplus\sp(1)$,  $\sp(r,\Real)\oplus\sl(2,\Real)$,
 $\sp(r,\Co)\oplus\sl(2,\Co)$,
$\spin(7)\subset\so(8)$, $\spin(4,3)\subset(4,4)$,
$\spin(7,\Co)\subset\so(8,8)$, $G_2\subset\so(7)$,
$G^*_{2(2)}\subset\so(4,3)$ and $G^\Co_2\subset\so(7,7)$. This
result, especially in the case of Riemannian manifolds, has a lot
of consequences and applications both in geometry and physics, see
 \cite{Al1,Besse,Bryant2,Joyce07} and the references therein.

Since theoretical physicists discovered supersymmetry,
supermanifolds began to play an important role both in mathematics
and physics \cite{DelMor,Leites,Manin,Var}. There is also an
interest to Riemannian supermanifolds, see e.g.
\cite{Asorey,Blaga,Cirilo,Cortes1,Go}. In particular,  Calabi-Yau
supermanifolds (i.e. Riemannian supermanifolds with the holonomy
algebras contained in $\su(p_0,q_0|p_1,q_1)$) were considered
recently in several physical papers, e.g. \cite{ABDS,RW05,Zhou}.

Holonomy groups and holonomy algebras for superconnections on supermanifolds are introduced recently in \cite{SupHol}. In
the present paper we generalize the result of M.~Berger for the case of Riemannian supermanifolds.

For a supersubalgebra $\g\subset\osp(p,q|2m)$ we define the vector
superspace $\R(\g)$ of curvature tensors of type $\g$ that
consists of linear maps from $\Lambda^2\Real^{p,q|2m}$ to $\g$
satisfying the Bianchi super identity. We call $\g$ {\it a Berger
superalgebra} if $\g$ is spanned by the images of the elements of
$\R(\g)$. The holonomy algebra of a Riemannian supermanifold is a
Berger superalgebra. Consequently, Berger superalgebras may be
considered as the candidates to the holonomy algebras. A Berger
superalgebra $\g\subset\osp(p,q|2m)$ is called {\it symmetric} if
the superspace $\R^\nabla(\g)$, consisting of linear maps  from
$\Real^{p,q|2m}$ to $\R(\g)$ satisfying the second Bianchi super
identity, is trivial. If such $\g$ is the holonomy algebra of a
Riemannian supermanifold, then this supermanifold is locally
symmetric. The classification of irreducible symmetric Berger
superalgebras can be deduced from \cite{Serganova}. Irreducible
Berger supersubalgebras  $\g\subset\osp(0|2m)$ are classified in
\cite{Odd}.  In the present paper we classify irreducible
non-symmetric Berger supersubalgebras $\g\subset\osp(p,q|2m)$
($p+q>0$) of the form
\begin{equation}\label{vidg}\g=(\oplus_i\g_i)\oplus\z,\end{equation}
where $\g_i$ are simple Lie superalgebras of classical type and
$\z$ is a trivial or one-dimensional center. The obtained list is
the following:
\begin{align*} &\osp(p,q|2m),\quad \osp(r|2k,\Co),\quad
\u(p_0,q_0|p_1,q_1), \quad \su(p_0,q_0|p_1,q_1),\\
&\hosp(r,s|k),\quad \hosp(r,s|k)\oplus\sp(1),\quad \osp^{sk}(2k|r,s)\oplus\sl(2,\Real),\\
&\osp^{sk}(2k|r,\Co)\oplus\sl(2,\Co).
\end{align*}
This list generalizes the list of irreducible holonomy algebras of
non-locally symmetric pseudo-Riemannian manifolds, in the same
time, we do not get the analogs of the important holonomy algebras
$G_2$ and $\spin(7)$. We discuss the geometries of the
supermanifolds with each of the obtained holonomy algebras.

Remark that in general a semi-simple Lie superalgebra $\g$ has the
form $\g=\oplus_i(\g_i\otimes\Lambda(n_i))$, where $\g_i$ is a
simple Lie superalgebra (either of classical or of Cartan type),
and $\Lambda(n_i)$ is the Grassmann superalgebra of $n_i$
variables \cite{Kac}. Moreover, there exist representations of
solvable Lie superalgebras in vector superspaces of dimensions
greater then one \cite{Kac}. Thus we consider only a part of
irreducible subalgebras $\g\subset\osp(p,q|2m)$. For irreducible
subalgebras $\g\subset\so(p,q)$ the property \eqref{vidg} holds
automatically, i.e. our assumption is quit natural.

In Section \ref{HpR} we recall some facts about holonomy algebras
of pseudo-Riemannian manifolds. In Section \ref{SecSupLie} and
\ref{SecSupmf} we provide some preliminaries on Lie superalgebras
and supermanifolds. In Section \ref{secMain}  we formulate the
main result of the paper. The remaining sections are dedicated to
the proof.  In Section \ref{secR(g)} we describe the superspaces
$\R(\g)$ for subalgebras $\g\subset\osp(p,q|2m)$. Let
$V=\Real^{p,q|2m}$. We show that any element $R\in\R(\g)_{\bar 0}$
satisfies $$\pr_{\so(p,q)}\circ R|_{\Lambda^2 V_{\bar
0}}\in\R(\pr_{\so(p,q)}\g_{\bar 0}),\quad \pr_{\sp(2m,\Real)}\circ
R|_{\Lambda^2 V_{\bar 1}}\in\bar \R(\pr_{\sp(2m,\Real)}\g_{\bar
0}),$$ where $\R(\pr_{\so(p,q)}\g_{\bar 0})$ is the space of
curvature tensors for the subalgebra $\pr_{\so(p,q)}\g_{\bar
0}\subset\so(p,q)$ and $\bar R(\pr_{\sp(2m,\Real)}\g_{\bar 0})$ is
a similar space for the subalgebra $\pr_{\sp(2m,\Real)}\g_{\bar
0}\subset\sp(2m,\Real)$, this space is defined in Section
\ref{SBa}. Next, any $R\in\R(\g)_{\bar 1}$ satisfies
$$\pr_{\so(p,q)}\circ R(\cdot,\xi)|_{V_{\bar
0}}\in\P_\eta(\pr_{\so(p,q)}\g_{\bar 0}),\quad
\pr_{\sp(2m,\Real)}\circ R(\cdot,x)|_{V_{\bar 1}}\in
P_\omega(\pr_{\sp(2m,\Real)}\g_{\bar 0}),$$ where $x\in V_{\bar
0}$ and $\xi\in V_{\bar 1}$ are fixed,
$\P_\eta(\pr_{\so(p,q)}\g_{\bar 0})$ is the space of weak
curvature tensors for the subalgebra $\pr_{\so(p,q)}\g_{\bar
0}\subset\so(p,q)$ and $\P_\omega(\pr_{\sp(2m,\Real)}\g_{\bar 0})$
is a similar space for the subalgebra $\pr_{\sp(2m,\Real)}\g_{\bar
0}\subset\sp(2m,\Real)$. We discuss these spaces in Sections
\ref{WBa} and \ref{SBa}. Elements of the space $\P_\eta(\so(n))$
appear as a part of the curvature tensor of a Lorentzian manifold
\cite{Leistner,onecomp}. These properties of the space $\R(\g)$
give strong conditions on the representation $\g_{\bar
0}\subset\so(p,q)\oplus\sp(2m,\Real)$ of the even part of $\g$. We
prove that under some assumption $\pr_{\so(p,q)}\g_{\bar
0}\subset\so(p,q)$ is the holonomy algebra of a pseudo-Riemannian
manifold and $\pr_{\sp(2m,\Real)}\g_{\bar 0}\subset\sp(2m,\Real)$
is in the list of possible reductive holonomy algebras of
Riemannian odd supermanifolds. These facts allow us to prove the
classification theorem. Many facts and tables that we use are
collected in \cite{Odd}.

\section{Holonomy algebras of pseudo-Riemannian manifolds}\label{HpR}

The theory of holonomy algebras of pseudo-Riemannian manifolds can
be found e.g. in \cite{Besse,Bryant2,Joyce07}. Here we collect
some facts that motivate the topic of this paper and will be used
below.

Let $(M,g)$ be a connected pseudo-Riemannian manifold of signature
$(p,q)$. The holonomy group of $(M,g)$ at a point $x\in M$ is a
Lie group that consists of the pseudo-orthogonal transformations
given by the parallel displacements along piece-wise smooth loops
at the point $x$, and it can be identified with a Lie subgroup of
the pseudo-orthogonal Lie group ${\rm O}(p,q)$. The corresponding
Lie subalgebra of $\so(p,q)$ is called the holonomy algebras. If
the manifold $M$ is simply connected, then the holonomy group is
connected and it is uniquely defined by the holonomy algebra.

A subalgebra $\g\subset\so(p,q)$ is called {\it
weakly-irreducible} if it does not preserve any proper
non-degenerate subspace of the pseudo-Euclidean space
$\Real^{p,q}$. By the Wu Theorem \cite{Wu} any pseudo-Riemannian
manifold can be decomposed (at least locally) in the product of a
flat pseudo-Riemannian manifold and of pseudo-Riemannian manifolds
with irreducible holonomy algebra. In particular, for any holonomy
algebra $\g\subset\so(V)$ there exists an orthogonal decomposition
$$\Real^{p,q}=V_0\oplus V_1\oplus\cdots\oplus V_r$$ into a direct sum of
pseudo-Euclidean subspaces and a decomposition
$$\g=\g_1\oplus\cdots\oplus \g_r$$ into a direct sum of ideals such that $\g_i$
annihilates $V_j$ if $i\neq j$ and $\g_i\subset \so(V_i)$ is the
weakly-irreducible holonomy algebra of a pseudo-Riemannian
manifold.

Possible connected irreducible holonomy groups (i.e. possible
irreducible holonomy algebras) of not locally symmetric
pseudo-Riemannian manifolds classified Marcel Berger in
\cite{Ber}.  Later it was proved that all these algebras  can be
realized as the holonomy algebras of pseudo-Riemannian manifolds
\cite{Bryant2}. Here is the list and the description of the
corresponding geometries (see e.g. \cite{Bryant2}):

\begin{itemize}
\item[] $\so(p,q)$: generic pseudo-Riemannian manifolds;
\item[] $\so(p,\Co)\subset\so(p,p)$: generic holomorphic pseudo-Riemannian manifolds;
\item[] $\u(r,s)\subset\so(2r,2s)$: pseudo-K\"ahlerian manifolds;
\item[] $\su(r,s)\subset\so(2r,2s)$: special pseudo-K\"ahlerian manifolds or Calabi-Yau manifolds;
\item[] $\sp(r,s)\subset\so(4r,4s)$: pseudo-hyper-K\"ahlerian manifolds;
\item[] $\sp(r,s)\oplus\sp(1)\subset\so(4r,4s)$: pseudo-quaternionic-K\"ahlerian manifolds;
\item[] $\sp(r,\Real)\oplus\sl(2,\Real)\subset\so(2r,2r)$: pseudo-para-K\"ahlerian manifolds;
\item[] $\sp(r,\Co)\oplus\sl(2,\Co)\subset\so(4r,4r)$: complex pseudo-para-K\"ahlerian manifolds;
\item[] $\spin(7)\subset\so(8)$, $\spin(4,3)\subset\so(4,4)$, $\spin(7,\Co)\subset\so(8,8)$:
8-dimensional pseudo-Riemannian manifolds with a parallel 4-form
and their complex analog;
\item[] $G_2\subset\so(7)$, $G^*_{2(2)}\subset\so(4,3)$,
$G^\Co_2\subset\so(7,7)$: 7-dimensional pseudo-Riemannian
manifolds with a parallel 3-form and their complex analog.
\end{itemize}

Compact Riemannian manifolds with the holonomy groups $\SU(2)$,
$\SU(3)$, $G_2$ and ${\rm Spin}(7)$ are extremely useful in
theoretical physics, see references in~\cite{Joyce07}.

The list of weakly-irreducible (and irreducible)  reductive
holonomy algebras of locally symmetric pseud-Riemannian manifolds
can be found in \cite{Ber57}.

In the case of Riemannian manifolds  weakly-irreducible
subalgebras $\g\subset\so(n)$ are the same as the irreducible
ones. Weakly-irreducible not irreducible Berger subalgebras
$\g\subset\so(p,q)$ are classified only if $p=1$  ($p$ is the
number of minuses of the metric), i.e. in the case of Lorentzian
manifolds, if $\g\subset\u(1,q)\subset\so(2,2q)$ and there are
some partial results in the neural signature $(p,p)$, see the
review \cite{IRMA}.

Let us shortly explain how a classification of holonomy algebras
can be obtained.

Let $\g\subset\so(p,q)=\so(V)$ be a subalgebra, where
$V=\Real^{p,q}$ is the pseudo-Euclidean space. The space of
curvature tensors of type $\g$ is defined as follows
$$\R(\g)=\left\{R\in\Lambda^2 V^*\otimes
\g\,\left|\begin{array}{c}R(X, Y)Z+R(Y, Z)X+R(Z, X)Y=0\\ \text{
}\text{ for all  } X,Y,Z\in V\end{array}\right\}\right..$$ The
above identity is called {\it the first Bianchi identity}. Let
$L(\R(\g))\subset\g$ be the ideal spanned by the images of the
elements form $\R(\g)$. From the Ambrose-Singer Theorem it follows
that for the holonomy algebra $\g$ of a pseudo-Riemannian manifold
it holds $L(\R(\g))=\g$.

Consider the vector space
$$\R^\na(\g)= \left\{S\in V^* \otimes\R(g)\left|\begin{matrix}S_X(Y,Z)+S_Y(Z,X)+S_Z(X,Y)=0\\
\text{for all  } X,Y,Z\in V
\end{matrix}\right\}\right..$$ If $\g$ satisfies $\R^\na(\g)=0$,
then any pseudo-Riemannian manifold with the holonomy algebra $\g$
is automatically locally flat. M.~Berger classified irreducible
subalgebras $\g\subset\so(p,q)$ that satisfy $L(\R(\g))=\g$ and
$\R^\na(\g)\neq 0$. By this reason subalgebra $\g\subset\so(p,q)$
 satisfying $L(\R(\g))=\g$ are called Berger algebras
 \cite{MS99,Sch}.
If $\g\subset\so(V)$ is a Berger subalgebra and $\R^\na(\g)=0$,
then $\g$ is called {\it a symmetric Berger algebra}, otherwise
$\g$ is called {\it a non-symmetric Berger algebra}. It is known
that if $\g$ is a reductive Lie algebra and the $\g$-module
$\R(\g)$ is trivial, then $\R^\na(\g)=0$ \cite{Sch}.

Let  $\g\subset\so(V)$ be a reductive weakly-irreducible Berger
subalgebra. If $\g$ is not irreducible, then it preserves a
degenerate subspace $W\subset V$. Consequently, $\g$ preserves the
isotropic subspace $L=W\cap W^\perp$. Since $\g$ is totally
reducible, there exists a complementary  invariant subspace
$L'\subset V$. Since $\g$ is weakly-irreducible, the subspace $L'$
is degenerate. If $L'$ is not isotropic, then $\g$ preserves the
kernel of the restriction of the metric to $L'$ and $\g$ preserves
a complementary subspace in $L'$ to this kernel, which is
non-degenerate. Hence $L'$ is isotropic and $V=L\oplus L'$ is the
direct sum of isotropic subspaces. This can happen only if $p=q$.
The metric on $V$ allows to identify $L'$ with the dual space
$L^*$ and the representations of $\g$ on $L$ and $L'$ are dual.
This shows that the representation $\g\subset\gl(L)$ is
irreducible. Let $R\in\R(\g)$. From the Bianchi identity it
follows that $R(x,y)=0$ and $R(\varphi,\psi)=0$ for all $x,y\in L$
and $\varphi,\psi\in L^*$. Moreover, for each fixed $\varphi\in
L^*$ it holds $R(\cdot,\varphi)\in(\g\subset\gl(L))^{(1)}$, where
$(\g\subset\gl(L))^{(1)}$ is the first prolongation for the
representation $\g\subset\gl(L)$ (similarly, for each fixed $x\in
L$ it holds $R(\cdot,x)\in(\g\subset\gl(L^*))^{(1)}$).
Consequently, $(\g\subset\gl(L))^{(1)}\neq 0$. The list of
irreducible subalgebras  and $\g\subset\gl(L)$ with
$(\g\subset\gl(L))^{(1)}\neq 0$ can be found in \cite{Bryant2}.

\section{Lie superalgebras}\label{SecSupLie}

Here we define some Lie superalgebras that we will use. Further
information about Lie superalgebras can be found e.g. in
\cite{DelMor,Kac,Leites}. For the purposes of this paper, in
\cite{Odd} we list simple Lie superalgebras of classical type.

{\it A vector superspace} $V$ is a $\mathbb{Z}_2$-graded vector
space $V=V_{\bar 0}\oplus V_{\bar 1}$, $\mathbb{Z}_2=\{\bar
0,\bar1\}$. The elements $x\in V_{\bar 0}\cup V_{\bar 1}$ are
called homogeneous. Elements $x\in V_{\bar 0}$ are called even,
and we define the parity $|x|$ of $x$ by putting $|x|=\bar 0$,
while elements $x\in V_{\bar 1}\backslash\{0\}$ are called odd, we
write $|x|=\bar 1$. All notions and constructions defined for the
usual vector spaces can be extended to the vector superspaces,
i.g. if $V$ and $W$ are vector superspaces, then their tensor
product $V\otimes W$ is a vector superspace with the grading
$$(V\otimes W)_{\bar 0}=(V_{\bar 0}\otimes W_{\bar
0})\oplus(V_{\bar 1}\otimes W_{\bar 1}),\quad (V\otimes W)_{\bar
1}=(V_{\bar 0}\otimes W_{\bar 1})\oplus(V_{\bar 1}\otimes W_{\bar
0}).$$ Linear maps $F:V\to W$ that preserve the parity are called
{\it morphisms}.

{\it A superalgebra} is a vector superspace $A=A_{\bar 0}\oplus
A_{\bar 1}$ with a morphism $\cdot:A\otimes A\to A$, i.e. for
homogeneous vectors  it holds $|x\cdot y|=|x|+|y|$. An important
example of a superalgebra is a {\it Grassmann superalgebra}
$$\Lambda(m)=\Lambda\Real^m=\oplus_{k=1}^m\Lambda^k\Real^m$$ with
the obvious $\mathbb{Z}_2$-grading. The Grassmann superalgebra is
super commutative, this means that its homogeneous elements
satisfy $xy=(-1)^{|x||y|}yx$.

{\it A Lie superalgebra} is a superalgebra $\g=\g_{\bar 0}\oplus
\g_{\bar 1}$ with the multiplication $[\cdot,\cdot]:\g\otimes
\g\to \g$ such that the homogeneous elements satisfy
$$[x,y]=-(-1)^{|x||y|}[y,x],$$ $$[[x,y],z]+(-1)^{|x|(|y|+|z|)}[[y,z],x]+(-1)^{|z|(|x|+|y|)}[[z,x],y]=0.$$
In particular, $\g_{\bar 0}$ is a Lie algebra and $\g_{\bar 1}$ is
a $\g_{\bar 0}$-module.

A simple Lie superalgebra $\g$ is of {\it classical type} if the
representation of $\g_{\bar 0}$ on $\g_{\bar 1}$ is totally
reducible, other wise it is of {\it the Cartan type}. A simple Lie
superalgebra $\g$ of classical type is of {\it type I} if the
representation of $\g_{\bar 0}$ in $\g_{\bar 1}$ is irreducible
and it is of {\it type II} if $\g_{\bar 1}$ is a direct sum of two
irreducible $\g_{\bar 0}$-modules.

Consider several examples.

Consider the vector superspace $V=V_{\bar 0}\oplus V_{\bar 1}$.
The Lie algebra of all endomorphisms of $V$ is denoted by
$\gl(V)$. Its even part $\gl(V)_{\bar 0}=\gl(V_{\bar
0})\oplus\gl(V _{\bar 1})$ consists of even endomorphisms, i.e. of
morphisms; the odd part of $\gl(V)$ consists of odd endomorphisms,
i.e. of the endomorphisms changing the parity and it holds
$\gl(V)_{\bar 1}=V^*_{\bar 0}\otimes V_{\bar 1}\oplus V^*_{\bar
1}\otimes V_{\bar 0}$.

Consider the vector superspace
$V=\Real^{n|m}=\Real^n\oplus\Pi\Real^m$. Here $\Pi$ is the parity
changing functor, in particular it is used to show that the vector
space $\Pi\Real^m$ is purely odd. In this case $\gl(V)$ is denote
by $\gl(n|m,\Real)$ (if $V$ is a real vector superspace). In the
matrix form we have
$$\gl(n|m,\Real)=\left\{\left(\begin{array}{cc}A&B\\
C&D\end{array} \right)\right\},$$ where $$\gl(n|m,\Real)_{\bar
0}=\left\{\left(\begin{array}{cc}A&0\\ 0&D\end{array}
\right)\right\},\quad \gl(n|m,\Real)_{\bar
1}=\left\{\left(\begin{array}{cc}0&B\\ C&0\end{array}
\right)\right\}.$$ Consider the Lie superalgebra
$\sl(n|m,\Real)=\{\xi\in\gl(n|m,\Real)|\str\xi=0\}$, where
$$\str\left(\begin{array}{cc}A&B\\ C&D\end{array}\right)=\tr A-\tr
D.$$ If $n\neq m$, then the Lie superalgebra $\sl(n|m)$ is simple
of the classical type II.

Let $g$ be an even  non-degenerate supersymmetric form on the
vector superspace $\Real^{n|2m}=\Real^n\oplus\Pi\Real^{2m}$, i.e.
$g(\Real^n,\Pi\Real^{2m})=g(\Pi\Real^{2m},\Real^n)=0$, the
restriction $\eta$ of $g$ to $\Real^n$ is non-degenerate and
symmetric (with some signature $(p,q),\,\,p+q=n$), and the
restriction $\omega$ of $g$ to $\Real^{2m}$ is non-degenerate and
skew-symmetric. Having such form, we denote $\Real^{n|2m}$ also by
$\Real^{p,q|2m}$. The orthosymplectic Lie superalgebra is defined
as the supersubalgebra of $\gl(n|2m,\Real)$ preserving $g$,
$$\osp(p,q|2m)_{\bar i}=\{\xi\in\gl(n|2m,\Real)_{\bar i}|\, g(\xi
x,y)+(-1)^{|x|\bar i}g(x,\xi y)=0 \},\,\,\bar i\in\mathbb{Z}_2.$$
If $p+q\neq 2$, then the Lie superalgebra $\osp(p,q|2m)$ is simple
of classical type I.
 In
particular, if the restriction of $g$ to $\Real^n$ is
positive definite, choose a basis with respect to which the matrix of $g$ has the form $\left(\begin{array}{ccc} 1_n&0&0\\
0&0&1_m\\ 0&-1_m&0\end{array}\right)$. Then, $$\osp(n|2m)=\left\{\left. \left(\begin{array}{ccc} A&B_1&B_2\\
B_2^t&C_1&C_2\\
-B_1^t&C_3&-C_1^t\end{array}\right)\right|A^t=-A,C_2^t=C_2,C_3^t=C_3\right\}.$$
Similarly, the Lie superalgebra $\osp^{sk}(2k|p,q)$, $p+q=m$, is
defined as the supersubalgebra of $\gl(2k|m,\Real)$  preserving an
even non-degenerate supersymmetric form on the vector superspace
$\Real^{2k}\oplus\Pi\Real^{m}$, in this case the restriction of
$g$ to $\Real^{2k}$ is non-degenerate and skew-symmetric, and the
restriction of $g$ to $\Real^m$ is non-degenerate and symmetric
(with the signature $(p,q),\,\,p+q=m$). For example,
$$\osp(2k|m)=\left\{\left.
\left(\begin{array}{ccc} C_1&C_2&B_1\\ C_3&-C_1^t&B_2\\
B_2^t&-B_1^t&A\end{array}\right)\right|A^t=-A,C_2^t=C_2,C_3^t=C_3\right\}.$$

Suppose that $n$ is even and suppose that we have an even
non-degenerate supersymmetric form on the vector superspace
$\Real^n\oplus\Pi\Real^{2m}$ such that the restriction of $g$ to
$\Real^n$ has signature $(2p_0,2q_0),\,\,2p_0+2q_0=n$. Suppose
that we have also a complex structure $J$ on
$\Real^n\oplus\Pi\Real^{2m}$ such that $g(Jx,Jy)=g(x,y)$ for all
$x,y\in\Real^n\oplus\Pi\Real^{2m}$. Note that the restriction of
$g(J\cdot,\cdot)$ to $\Real^{2m}$ is symmetric and non-degenerate
and let this form have the signature $(2p_1,2q_1),\,\,p_1+q_1=m$.
By definition,
$$\u(p_0,q_0|p_1,q_1)=\{\xi\in\osp(2p_0,2q_0|2m)|[\xi,J]=0\},$$
$$\su(p_0,q_0|p_1,q_1)=\{\xi\in\u(p_0,q_0|p_1,q_1)|\str(J\circ\xi)=0\}.$$
Similarly, suppose that $n$ and $2m$ are integers divided by 4,
$m=2k$, and suppose that we have an even  non-degenerate
supersymmetric form on the vector superspace
$\Real^n\oplus\Pi\Real^{2m}$ such that the restriction of $g$ to
$\Real^n$ has the signature $(4r,4s),\,\,4r+4s=n$. Suppose that we
have a quaternionic structure $I,J,K$ on
$\Real^n\oplus\Pi\Real^{2m}$ (i.e. $I,J,K$ are complex structures
and they generate the Lie algebra isomorphic to $\sp(1)$) such
that $g(Ix,Iy)=g(Jx,Jy)=g(Kx,Ky)=g(x,y)$ for all
$x,y\in\Real^n\oplus\Pi\Real^{2m}$.  By definition,
$$\hosp(r,s|k)=\{\xi\in\osp(4r,4s|2m)|[\xi,I]=[\xi,J]=[\xi,K]=0\}.$$
Note that the normalizer of $\sp(1)$ in $\osp(4r,4s|2m)$ coincides
with $\sp(1)\oplus \hosp(r,s|k)$. Remark also that
$\hosp(r,s|k)_{\bar 0}\cap\sp(2m,\Real)=\so^*(k)$, where
$\so^*(k)$ is the subalgebra of $\gl(k,\mathbb{H})$ acting on
$\Real^{4k}=\mathbb{H}^k$ and preserving the
skew-quaternionic-Hermitian form
$$\omega(\cdot,\cdot)+i\omega(I\cdot,\cdot)+j\omega(J\cdot,\cdot)+k\omega(K\cdot,\cdot),$$
here $\omega$ is the restriction of $g$ to $\Real^{4k}$ considered
as a skew-symmetric bilinear form on $\Real^{4k}$.

\section{Supermanifolds and their holonomy groups}\label{SecSupmf}

The introduction to the theory  of supermanifolds can be found in
\cite{DelMor,Leites,Manin,Var}. Riemannian supermanifolds are
considered e.g. in \cite{Asorey,Blaga,Cirilo,Cortes1,Go}. The
holonomy theory is introduced in \cite{SupHol}.

{\it  A real smooth  supermanifold} $\M$ of dimension $n|m$ is a
pair $(M,\O_\M)$, where $M$ is a Hausdorff topological space and
$\O_\M$ is a sheaf of commutative superalgebras with unity over
$\Real$ locally isomorphic to
$$\Real^{n;m}=(\Real^n,\O_{\Real^{n;m}}=\O_{\Real^n}\otimes
\Lambda(m)),$$ where $\O_{\Real^n}$ is the sheaf of smooth
functions on $\Real^n$ and $\Lambda(m)$ is the Grassmann
superalgebra of $m$ generators. The sections of the sheaf $\O_\M$
are called {\it superfunctions (or just functions)} on $\M$. The
ideal $(\O_\M)_{\bar 1}\oplus((\O_\M)_{\bar 1})^2$ consists of the
nilpotent elements of $\O_\M$, and the sheaf $\O_M$ defined as the
quotient $\O_\M/((\O_\M)_{\bar 1}\oplus((\O_\M)_{\bar 1})^2)$
furnish $M$ with the structure of a real  smooth manifold. We get
the canonical projection $\ti:\O_\M\to \O_M$, $f\mapsto \t f$.
{\it The value of a superfunction} $f$ at a point $x\in M$ is by
definition $\t f(x)$. If $m=0$, then $\M=M$ is a smooth manifold.

Denote by $\T_\M$ the tangent sheaf, i.e.  the sheaf of
superderivatives of the sheaf $\O_\M$. For a point $x\in M$ define
the tangent space $$T_x\M=\T_\M(U)/(\O_\M(U)_x\T_\M(U)),$$ where
$U\subset M$ is any open subset containing $x$ and $\O_\M(U)_x$ is
the ideal in $\O_\M(U)$ consisting of functions vanishing at the
point $x$. The space $T_x\M$ is a vector superspace of the same
dimension as $\M$ and its even part coincides with $T_xM$.
Consider the tangent bundle $T\M=\cup_{x\in M}T_x\M$ as a bundle
of vector superspaces over $M$, its even part is the tangent
bundle of $M$. We may consider also $T\M$ as a usual vector bundle
over $M$. For a section $X\in\T_\M(U)$, its value $X_x\in T_x\M$
is defined in the obvious way. We get the natural projection
$\ti:\T_\M(U)\to\Gamma(T\M,U)$.

{\it A connection} on $\M$ is an even morphism
$\na:\T_\M\otimes_\Real \T_\M\to\T_\M$ of sheaves of supermodules
over $\Real$ such that
$$\na_{fY}X=f\na_YX\quad\text{and}\quad\na_YfX=(Yf)X+(-1)^{|Y||f|}f\na_YX$$
for all homogeneous functions $f$, vector fields $X,Y$ on $\M$,
here $|\cdot|\in\mathbb{Z}_2=\{{\bar 0},{\bar 1}\}$ denotes the
parity. In particular, $|\na_YX|=|Y|+|X|$. The {\it curvature
tensor} of the connection $\na$ is given by
\begin{equation}\label{R} R(Y,Z)=[\na_Y,\na_Z]-\na_{[Y,Z]},\end{equation}
where $Y$ and $Z$ are vector fields on $\M$.

Obviously, the restriction $$\t\na=(\na|_{\Ga(TM)\otimes
\Ga(T\M)})^\ti:\Ga(TM)\otimes \Ga(T\M)\to\Ga(T\M)$$ is a
connection on the vector bundle  $T\M$. Since $\na$ is even, the
subbundles $TM=(T\M)_{\bar 0},(T\M)_{\bar 1}\subset T\M$ are
parallel. We obtain also a linear connection on $M$.
 Let $\ga:[a,b]\subset\Real\to M$ be a curve and $\tau_\gamma:T_{\ga(a)}\M\to T_{\ga(b)}\M$
 the parallel displacement
along $\ga$ defined by the connection $\t\nabla$.

{\it  The holonomy algebra} $\hol(\na)_x$ of the connection $\na$
at a point $x\in M$ is the supersubalgebra of the Lie superalgebra
$\gl(T_x\M)$ generated by the operators of the form
$$\tau_\ga^{-1}\circ
\na^r_{{Y_r},...,{Y_1}}R_y(Y,Z)\circ\tau_\ga:T_x\M\to T_x\M,$$
where $\gamma$ is any piecewise smooth curve in $M$ beginning at
the point $x$; $y\in M$ is the end-point of the curve $\gamma$,
$r\geq 0$, and  $Y,Z,Y_1,...,Y_r\in T_y\M$.

Now we define the holonomy group.  Recall that a {\it Lie
supergroup} $\mathcal G=(G,\O_\mathcal G)$ is a group object in
the category of supermanifolds. The underlying smooth manifold $G$
is a Lie group. The Lie superalgebra $\g$ of $\mathcal G$ can be
identified with the tangent space to $\mathcal G$ at the identity
$e\in G$. The Lie algebra of the Lie group $G$ is the even part
$\g_{\bar 0}$ of the Lie superalgebra $\g$. Any Lie supergroup
$\mathcal G$ is uniquely given by a pair $(G,\g)$ ({\it
Harish-Chandra pair}), where $G$ is a Lie group, $\g=\g_{\bar
0}\oplus\g_{\bar 1}$ is a Lie superalgebra such that $\g_{\bar 0}$
is the Lie algebra of the Lie group $G$ and a certain conditions
hold  \cite{DelMor,Go}.

Denote by $\Hol(\na)^0_x$ the connected Lie subgroup of the Lie
group $\GL((T_x\M)_{\bar 0})\times \GL((T_x\M)_{\bar 1})$
corresponding to the Lie subalgebra $(\hol(\na)_x)_{\bar 0}\subset
\gl((T_x\M)_{\bar 0})\oplus \gl((T_x\M)_{\bar
1})\subset\gl(T_x\M)$. Let $\Hol(\na)_x$ be the Lie subgroup of
the Lie group $\GL((T_x\M)_{\bar 0})\times \GL((T_x\M)_{\bar 1})$
generated by the Lie groups $\Hol(\na)^0_x$ and $\Hol(\t\na)_x$.
Clearly, the Lie algebra of the Lie group $\Hol(\na)_x$ is
$(\hol(\na)_x)_{\bar 0}$. The Lie supergroup $\Hols(\na)_x$ given
by the Harish-Chandra pair $(\Hol(\na)_x,\hol(\na)_x)$ is called
the holonomy group of the connection $\na$ at the point $x$. The
Lie supergroup $\Hols(\na)^0_x$ given by the Harish-Chandra pair
$(\Hol(\na)^0_x,\hol(\na)_x)$ is called the restricted holonomy
group of the connection $\na$ at the point $x$.

If $\M=M$ is a usual smooth manifold and $\nabla$ is a connection
on $M$, than $\Hols(\na)_x$ and $\hol(\na)_x$ coincide with the
holonomy group and the holonomy algebra of $(M,\nabla)$,
respectively.

If the manifold $M$ is simply connected, then the group
$\Hol(\na)_x$ is connected and the whole information about
$\Hols(\na)_x$ gives us the holonomy algebra $\hol(\na)_x$.

{\it A Riemannian supermanifold} $(\M,g)$ is a supermanifold $\M$
of dimension $n|2m$ endowed with an even non-degenerate
supersymmetric metric $$g:\T_\M\otimes_{\O_\M}\T_\M\to\O_\M,$$ see
 e.g. \cite{Cortes1}.  In particular, the value $g_x$ of $g$ at a
point $x\in M$ satisfies: $g_x((T_x\M)_{\bar 0},(T_x\M)_{\bar
1})=0$, $g_x|_{(T_x\M)_{\bar 0}\times (T_x\M)_{\bar 0}}$ is
non-degenerate, symmetric and $g_x|_{(T_x\M)_{\bar 1}\times
(T_x\M)_{\bar 1}}$ is non-degenerate, skew-symmetric. The metric
$g$ defines a pseudo-Riemannian metric $\tilde g$ on the manifold
$M$. Note that $\tilde g$ is not assumed to be positively defined.
The supermanifold $(\M,g)$ has a unique linear connection  $\na$
such that $\na$ is torsion-free and $\na g=0$. This connection is
called {\it the Levi-Civita connection}. We denote the holonomy
algebra of the connection $\na$ by $\hol(\M,g)$. Since $g$ is
parallel, $\hol(\M,g)\subset\osp(p,q|2m)$ and
$\Hol(\tilde\na)\subset \Ort(p,q)\times \Sp(2m,\Real)$, where
$(p,q)$ is the signature of the pseudo-Riemannian metric $\tilde
g$.

We call a supersubalgebra $\g\subset \osp(p,q|2m)$ {\it
weakly-irreducible} if it does not preserve any non-degenerate
vector supersubspace of $\Real^{p,q}\oplus\Pi\Real^{2m}$. The
following theorem generalizes the Wu theorem \cite{Wu}.

\begin{theorem}\cite{SupHol}\label{ThWuSuper} Let $(\M,g)$ be a Riemannian supermanifold such that the pseudo-Riemannian manifold $(M,\tilde g)$ is simply connected
and geodesically complete. Then there exist Riemannian
supermanifolds $(\M_0,g_0),(\M_1,g_1),...,(\M_r,g_r)$ such that
\begin{equation}\label{Mdec} (\M,g)=(\M_0\times\M_1\times\cdots\times\M_r,g_0+g_1+\cdots+g_r),\end{equation}
the supermanifold $(\M_0,g_0)$ is flat and the holonomy algebras
of the supermanifolds $(\M_1,g_1)$,...,$(\M_r,g_r)$ are
weakly-irreducible. In particular,
\begin{equation}\hol(\M,g)=\hol(\M_1,g_1)\oplus\cdots\oplus\hol(\M_r,g_r).\end{equation}

For general $(\M,g)$ decomposition \eqref{Mdec} holds locally.
\end{theorem}

\section{Berger superalgebras}\label{secMain}

Let $V=V_{\bar 0}\oplus V_{\bar 1}$ be a real or complex vector
superspace and  $\g\subset\gl(V)$ a supersubalgebra. {\it The
space of algebraic curvature tensors of type} $\g$ is the vector
superspace $\R(\g)=\R(\g)_{\bar 0}\oplus\R(\g)_{\bar 1},$ where
$$\R(\g)=\left\{R\in \Lambda^2 V^*
\otimes\g\left|\begin{matrix}R(X,Y)Z+(-1)^{|X|(|Y|+|Z|)}R(Y,Z)X\\+(-1)^{|Z|(|X|+|Y|)}R(Z,X)Y=0\\\text{for
all homogeneous } X,Y,Z\in V
\end{matrix}\right\}\right..$$ Here $|\cdot|\in\mathbb{Z}_2$ denotes the parity. The identity that satisfy the elements $R\in\R(\g)$
 is called {\it the first Bianchi super identity}.
Obviously, $\R(\g)$ is a $\g$-module with respect to the action
$A\cdot R=R_A$,
\begin{equation}\label{R_A}
R_A(X,Y)=[A,R(X,Y)]-(-1)^{|A||R|}R(AX,Y)-(-1)^{|A|(|R|+|X|)}R(X,AY),\end{equation}
where $A\in\g$, $R\in\R(\g)$ and $X,Y\in V$ are homogeneous.

If $\M$ is a supermanifold and $\na$ is a linear torsion-free
connection on the tangent sheaf $\T_\M$ with the holonomy algebra
$\hol(\na)_x$ at some point $x$, then for the covariant
derivatives of the curvature tensor  we have
$(\na^r_{Y_r,...,Y_1}R)_x\in\R(\hol(\na)_x)$ for all $r\geq 0$ and
tangent vectors  $Y_1,...,Y_r\in T_xM$. Moreover,
$|(\na^r_{Y_r,...,Y_1}R)_x|=|Y_1|+\cdots+|Y_r|$, whenever
$Y_1,...,Y_r$ are homogeneous.

Define the vector supersubspace $$L(\R(\g))=\spa\{R(X,Y)|R\in\R(\g),\,\,X,Y\in V\}\subset \g.$$ From \eqref{R_A} it
follows that $L(\R(\g))$ is an ideal in $\g$. We call a supersubalgebra $\g\subset\gl(V)$ {\it a Berger superalgebra} if
$L(\R(\g))=\g$.

If $V$ is a vector space, which can be considered as a vector superspace with the trivial odd part,  then
$\g\subset\gl(V)$ is a usual Lie algebra, which can be considered as a Lie superalgebra with the trivial odd part. Berger
superalgebras in this case are the same as the usual Berger algebras.

\begin{prop}  \cite{SupHol} Let $\M$ be a supermanifold of dimension $n|m$ with a linear torsion-free connection $\na$.
Then its holonomy algebra $\hol(\na)\subset\gl(n|m,\Real)$ is a Berger superalgebra.\end{prop}

Consider the vector superspace $$\R^\na(\g)= \left\{S\in V^*
\otimes\R(g)\left|\begin{matrix}S_X(Y,Z)+(-1)^{|X|(|Y|+|Z|)}S_Y(Z,X)\\+(-1)^{|Z|(|X|+|Y|)}S_Z(X,Y)=0\\
\text{for all homogeneous } X,Y,Z\in V
\end{matrix}\right\}\right..$$

If $\M$ is a supermanifold and $\na$ is a linear torsion-free
connection on $\T_\M$, then
$(\na^r_{Y_r,...,Y_2,\cdot}R)_x\in\R^\na(\hol(\na)_x)$ for all
$r\geq 1$ and $Y_2,...,Y_r\in T_xM$. Moreover,
$|(\na^r_{Y_r,...,Y_2,\cdot}R)_x|=|Y_2|+\cdots+|Y_r|$, whenever
$Y_2,...,Y_r$ are homogeneous.

A Berger superalgebra $\g$ is called {\it symmetric} if $\R^\na(\g)=0$. This is a generalization of the usual symmetric
Berger algebras, see e.g. \cite{Sch}, and the following is a generalization of the well-known fact about smooth manifolds.

\begin{prop}\cite{SupHol} Let $\M$ be a  supermanifold with a torsion free connection $\na$. If $\hol(\na)$ is a symmetric Berger
superalgebra, then $(\M,\na)$ is locally symmetric (i.e. $\nabla R=0$). If $(\M,\na)$ is a locally symmetric superspace,
then its curvature tensor at any point is annihilated by the holonomy algebra at this point and its image coincides with
the holonomy algebra.
\end{prop}

A geometric theory of Riemannian symmetric superspaces is
developed recently in \cite{Go}.

The proof of the following proposition is as in \cite{Sch}.

\begin{prop} Let $\g\subset\gl(V)$ be an irreducible Berger superalgebra of the form \eqref{vidg}. If $\g$ annihilates the module $\R(\g)$,
then $\g$ is a symmetric Berger superalgebra.
\end{prop}

In \cite{Serganova} simply connected symmetric superspaces of
simple Lie supergroups of isometries are classified. In particular
this implies the classification of the holonomy algebras of
Riemannian symmetric superspaces and of irreducible Berger
superalgebras $\g\subset \osp(p,q|2m)$ of the form \eqref{vidg}.
Hence we assume that the Riemannian supermanifolds under the
consideration are not locally symmetric.

\section{The Main Results}\label{secMain}

Here is the Main Theorem of this paper.

\begin{theorem}\label{SBTh} Let $(\M,g)$ be a not locally symmetric Riemannian supermanifold
of dimension $p+q|2m$ ($p+q>0$) with an irreducible holonomy
algebra
 $\hol(\M,g)\subset \osp(p,q|2m)$ that is of the form \eqref{vidg}, then
$\hol(\M,g)\subset\osp(p,q|2m)$ coincides with one of the Lie
superalgebras from Table~\ref{tabSuperB}.
\end{theorem}

\begin{tab}\label{tabSuperB}
Irreducible non-symmetric Berger supersubalgebras $\g\subset
\osp(p,q|2m)$ ($p+q>0$) of the from \eqref{vidg} and the connected
Lie subgroups $G\subset\SO(p,q)\times \Sp(2m,\Real)$ corresponding
to $\g_{\bar 0}\subset\so(p,q)\oplus\sp(2m,\Real)$.
$$\begin{array}{|c|c|c|} \hline \g&G&(p,q|2m)\\ \hline
\osp(p,q|2m)&\SO(p,q)\times \Sp(2m,\Real)&(p,q|2m)\\
 \osp(p|2k,\Co)&\SO(p,\Co)\times \Sp(2k,\Real)& (p,p|4k)\\
\u(p_0,q_0|p_1,q_1)&\Un(p_0,q_0)\times \Un(p_1,q_1)& (2p_0,2q_0|2p_1+2q_1)\\
\su(p_0,q_0|p_1,q_1)&\Un(1)(\SU(p_0,q_0)\times \SU(p_1,q_1))& (2p_0,2q_0|2p_1+2q_1)\\
\hosp(r,s|k)&\Sp(p_0,q_0)\times \SO^*(k)&(4r,4s|4k)\\
\hosp(r,s|k)\oplus\sp(1)&\Sp(1)(\Sp(p_0,q_0)\times \SO^*(k))&(4r,4s|4k)\\
\osp^{sk}(2k|r,s)\oplus\sl(2,\Real)&\Sp(2k,\Real)\times\SO(r,s)\times\SL(2,\Real) &(2k,2k|2r+2s)\\
\osp^{sk}(2k|r,\Co)\oplus\sl(2,\Co)&\Sp(2k,\Co)\times\SO(r,\Co)\times\SL(2,\Co)&(4k,4k|4r)\\
\hline\end{array}$$
\end{tab}

The Ricci tensor of a supermanifold $(\M,g)$ is defined by the
formula \begin{equation}\label{defRic}\Ric(Y,Z)=\str(X\mapsto
(-1)^{|X||Z|}R(Y,X)Z),\end{equation} where $U\subset M$ is open
and $X,Y,Z\in\T_\M(U)$ are homogeneous.

The definitions of the supermanifolds considered below are similar
to the usual ones, see e.g. \cite{Cortes1,SupHol}. Foe example, a
Riemannian supermanifold $(\M,g)$ is called a {\it K\"ahlerian
supermanifold} if it admits an even parallel $g$-orthogonal
complex structure. Since the holonomy algebra annihilates the
values of the parallel tensors \cite{SupHol}, in this case it must
be contained in $\u(p_0,q_0|p_1,q_1)$. By definition, {\it a
special K\"ahlerian supermanifold or a Calabi-Yau supermanifold}
is a Ricci-flat K\"ahlerian supermanifold.

\begin{prop}  \cite{SupHol} Let $(\M,g)$ be a  K\"ahlerian supermanifold, then $\Ric=0$
if and only if $\hol(\M,g)\subset\su(p_0,q_0|p_1,q_1)$. \end{prop}

Riemannian supermanifolds with the holonomy algebras
$\osp(p,q|2m)$ are generic.  Here we give the geometric
characterization of simply connected supermanifolds with the
holonomy algebras $\g$ different from $\osp(p,q|2m)$ (the
conditions on the holonomy algebra and the corresponding
geometries are equivalent for simply connected supermanifolds):

\begin{itemize}
\item[]$\g\subset\osp(p|2k,\Co)$: holomorphic Riemannian supermanifolds;
\item[]$\g\subset\u(p_0,q_0|p_1,q_1)$: K\"ahlerian supermanifolds;
\item[]$\g\subset\su(p_0,q_0|p_1,q_1)$: special K\"ahlerian supermanifolds or
Calabi-Yau supermanifolds;
\item[]$\g\subset\hosp(r,s|k)$: hyper-K\"ahlerian supermanifolds;
\item[]$\g\subset\hosp(r,s|k)\oplus\sp(1)$: quaternionic-K\"ahlerian supermanifolds;
\item[]$\g\subset\osp^{sk}(2k|r,s)\oplus\sl(2,\Real)$: para-K\"ahlerian
supermanifolds;
\item[]$\g\subset\osp^{sk}(2k|r,\Co)\oplus\sl(2,\Co)$: holomorphic para-K\"ahlerian
supermanifolds.
\end{itemize}

\begin{prop}  \cite{SupHol}
 Let $(\M,g)$ be a  quaternionic-K\"ahlerian supermanifold, then $\Ric=0$ if and only if
$\hol(\M,g)\subset\hosp(p_0,q_0|p_1,q_1)$. In particular, if
$(\M,g)$ is hyper-K\"ahlerian, then $\Ric=0$; if $M$ is simply
connected, $(\M,g)$ is quaternionic-K\"ahlerian and $\Ric=0$, then
$(\M,g)$ is hyper-K\"ahlerian. \end{prop}

Now the natural problem is to construct examples of supermanifolds
with each of the obtained possible holonomy algebras. Note that
examples of special K\"ahlerian manifolds (i.e. Calabi-Yau
manifolds) delivers the Calabi-Yau Theorem. In \cite{RW05} it is
shown that the Calabi-Yau Theorem does not hold for K\"ahlerian
supermanifolds of real odd dimension two. In \cite{Zhou} it is
shown that the arguments of \cite{RW05} work only for the odd
dimension two and it is conjectured that Calabi-Yau Theorem is
true for manifolds of odd dimensions bigger then two. In
\cite{ABDS} some examples of Calabi-Yau supermanifolds are
constructed.  Examples of quaternionic-K\"ahlerian supermanifolds
are constructed in \cite{Cortes1}.

\vskip0.2cm

The rest sections are dedicated to the proof of Theorem
\ref{SBTh}.

\section{Weak-Berger algebras}\label{WBa}

 Let $\g\subset\so(p,q)=\so(V)$ be a subalgebra. Denote by $\eta$ the pseudo-Euclidian metric on $V=\Real^{p,q}$.
 The vector space
$$\P_\eta(\g)=\left\{P\in V^*\otimes \g\,\left|\,\begin{array}{c} \eta(P(X)Y,Z)+\eta(P(Y)Z,X)+\eta(P(Z)X,Y)=0\\ \text{ for
all }X,Y,Z\in V\end{array}\right\}\right.$$ is called {\it the
space of  weak-curvature tensors of type } $\g$. A subalgebra
$\g\subset\so(V)$ is called a weak-Berger algebra if  $\h$ is
spanned by the images of the elements $P\in\P(\h)$. It is not hard
to see that if $R\in\R(\g)$ and $x\in V$ is fixed, then
$R(\cdot,x)\in\P_\eta(\g)$. In particular, any Berger algebra is a
weak-Berger algebra. The converse statement is not obvious, it is
proved recently in \cite{Leistner}.

\begin{theorem}\label{Leistner} Let $\g\subset\so(p,q)$ be an irreducible weak-Berger subalgebra, then it is a Berger subalgebra.
\end{theorem}

Remark that in the origin theorem $\g$ is a subalgebra of
$\so(n)$. The above result immediately follows from  the
complexification process described in \cite{Leistner}.

The spaces $\P_\eta(\g)$ for irreducible Berger subalgebras
$\g\subset\so(n)$ are found in \cite{onecomp}. This result can be
easily extended to the case of subalgebras $\g\subset\so(p,q)$. In
particular, it is proved that if
$\g=\g^1\oplus\g^2\subset\so(V^1\otimes V^2)=\so(V)$, where
$\g^1\subset\gl(V^1)$ and $\g^2\subset\gl(V^2)$ are irreducible,
then $$\P_\eta(\g^1\subset\so(V))=\P_\eta(\g^2\subset\so(V))=0,$$
unless the complexification of $\g\subset\so(V)$ coincides with
$$\sp(2m,\Co)\oplus\sl(2,\Co)\subset\so(4m,\Co).$$ For example, for
$\sp(m)\oplus\sp(1)\subset\so(4m)$ it holds
$$\P_\eta(\sp(1)\subset\so(V))=0,\quad
\P_\eta(\sp(m)\subset\so(V))=(\sp(2m,\Co)\subset\sl(2m,\Co))^{(1)}.$$

In  \cite{onecomp} it is shown that if $\g\subset\so(n)$ is an
irreducible subalgebra and $\P_\eta(\g)\neq 0$ or $\R(\g)\neq 0$,
then either $\g$ is a Berger subalgebra, or
$\g=\sp(\frac{n}{4})\oplus\Real J$. A similar result holds for
irreducible subalgebras $\g\subset\so(p,q)$.

Let $\g\subset\gl(L)$ be an irreducible subalgebra. Then $\g$ is a
weakly-irreducible subalgebra of $\so(L\oplus L^*)=\so(p,p)$,
where $p=\dim L$. Let $\eta$ be the natural metric on $L\oplus
L^*$. It is easy to see that $P\in\P_\eta(\g)$ if and only if
$$\pr_{\gl(L)}\circ P|_{L}\in (\g\subset\gl(L))^{(1)},\quad
\pr_{\gl(L^*)}\circ P|_{L}\in(\g\subset\gl(L^*))^{(1)}.$$ Thus if
$\g\subset\so(L\oplus L^*)$ is a weak-Berger subalgebra, then
$(\g\subset\gl(L))^{(1)}\neq \{0\}$ and all $\g\subset\gl(L)$ are
known \cite{Bryant2}.

\section{The case of Riemannian odd supermanifolds}\label{SBa}

Let $(\M,g)$ be a Riemannian supermanifold of dimension $0|2m$,
such supermanifolds are called odd. In this case for the holonomy
algebra we get $$\g\subset\osp(0|2m)\simeq\sp(2m,\Real),$$ i.e.
$\g$ is a usual Lie algebra acting in a purely odd vector
superspace. The possible irreducible holonomy algebras of such
supermanifolds are classified in \cite{Odd}.

Let $V$ be a real or complex vector space and $\g\subset\gl(V)$ a
subalgebra. The space of skew-symmetric curvature tensors of type
$\g$ is defined as follows
$$\bar\R(\g)=\left\{R\in\odot^2 V^*\otimes
\g\,\left|\begin{array}{c}R(X, Y)Z+R(Y, Z)X+R(Z, X)Y=0\\ \text{
}\text{ for all  } X,Y,Z\in V\end{array}\right\}\right..$$ The
subalgebra $\g\subset\gl(V)$ is called {\it a skew-Berger
subalgebra} if it is spanned by the images of the elements
$R\in\bar\R(\g)$. Obviously $\bar\R(\g)=\R(\g\subset\gl(\Pi V))$
and $\g\subset\gl(V)$ is a skew-Berger algebra if and only if
$\g\subset\gl(\Pi V)$ is a Berger superalgebra.

Let $\omega$ be the standard symplectic form on $V=\Real^{2m}$. A
subalgebra $\g\subset\sp(2m,\Real)$ is called {\it
weakly-irreducible} if it does not preserve any proper
non-degenerate subspace of $\Real^{2m}$. The Wu Theorem for
supermanifolds implies the following statement.
  Let $\g\subset\sp(2m,\Real)$ be an irreducible
skew-Berger subalgebra, then there is a decomposition
$$V=V_0\oplus V_1\oplus\cdots\oplus V_r$$ into a direct sum of
symplectic subspaces and a decomposition
$$\g=\g_1\oplus\cdots\oplus \g_r$$ into a direct sum of ideals
such that $\g_i$ annihilates $V_j$ if $i\neq j$ and $\g_i\subset
\sp(V_i)$ is a weakly-irreducible skew-Berger subalgebra.

Irreducible skew-Berger subalgebras $\g\subset\gl(n,\Co)$ are
classified  in \cite{Skew-Berger}. Irreducible skew-Berger
subalgebras $\g\subset\sp(2m,\Real)$ are classified in \cite{Odd}.

Let  $\g\subset\sp(2m,\Real)=\sp(V)$ be a reductive
weakly-irreducible subalgebra. Suppose that $\g$ is not
irreducible. As in Section \ref{HpR} we may show that $V$ is of
the form $V=L\oplus L^*$, where $\g\subset\gl(L)$ is irreducible.
If $\g\subset\sp(V)$ is a skew-Berger subalgebra, then $(\g\subset
\gl(L))^{[1]}\neq\{0\}$, where
$$\g^{[1]}=\{\varphi\in L^*\otimes\g|\varphi(x)y=-\varphi(y)x\text{
for all } x,y\in L\}$$ is the skew-symmetric prolongation of the
subalgebra $\g\subset\gl(L)$. Irreducible subalgebras
$\g\subset\gl(L)$ with $\g^{[1]}\neq 0$ are classified in
\cite{Odd}.

\vskip0.5cm

Let $V$ be a complex or real vector space with a symplectic form
$\omega$.
 Let $\g\subset\sp(V)$ be a subalgebra.
 The vector space
$$\P_\omega(\g)=\left\{P\in V^*\otimes \g\,\left|\,\begin{array}{c} \omega(P(X)Y,Z)+\omega(P(Y)Z,X)+\omega(P(Z)X,Y)=0\\ \text{
for all }X,Y,Z\in V\end{array}\right\}\right.$$ is called {\it the
space of  skew-symmetric weak-curvature tensors of type } $\g$. A
subalgebra $\g\subset\sp(V)$ is called a skew-symmetric
weak-Berger algebra if  $\g$ is spanned by the images of the
elements $P\in\P_\omega(\g)$. It is not hard to see that if
$R\in\bar\R(\g)$ and $X\in V$ is fixed, then
$R(\cdot,X)\in\P_\omega(\g)$. In particular, any skew-Berger
algebra is a skew-symmetric weak-Berger algebra. The converse
statement gives the following theorem.

\begin{theorem} Let $\g\subset\sp(2m,\Real)$ be an irreducible skew-symmetric weak-Berger subalgebra,
then it is a skew-Berger subalgebra.
\end{theorem}

{\it The proof} of this theorem is a modified copy of the proof
from \cite{Leistner} of Theorem \ref{Leistner}.

If the representation $\g\subset\sp(2m,\Real)$ is not absolutely
irreducible, then $\P_\omega(\g)$ is isomorphic to
$(\g_\Co\subset\gl(m,\Co))^{[1]}$ and the proof follows from
\cite{Skew-Berger,Odd}.

If the representation $\g\subset\sp(2m,\Real)$ is absolutely
irreducible, then we need a classification of irreducible
skew-symmetric weak-Berger subalgebras $\g\subset\sp(2m,\Co)$. It
can be achieve in the same way as the classification of
irreducible weak-Berger subalgebras $\g\subset\so(n,\Co)$. In
fact, in \cite{Leistner} it is obtained a necessary condition for
an irreducible subalgebra $\g\subset\so(n,\Co)$ to be a
weak-Berger subalgebra, then there were classified all subalgebras
$\g\subset\gl(n,\Co)$ satisfying this condition and the
subalgebras $\g\subset\sp(2m,\Co)$ were noted. It is easy to see
that a skew-symmetric weak-Berger subalgebras
$\g\subset\sp(2m,\Co)$ satisfy the same necessary condition. Thus
the proof follows immediately. $\Box$

\vskip0.5cm

The spaces $\P_\omega(\g)$ for irreducible weak-Berger subalgebras
$\g\subset\sp(2n,\Real)$ can be found by methods of
\cite{onecomp}.
 In particular, it is can be  proved that
 if $$\g=\g^1\oplus\g^2\subset\sp(V^1\otimes V^2)=\sp(V),$$ where $\g^1\subset\gl(V^1)$
  and $\g^2\subset\gl(V^2)$ are irreducible,
 then
 $$\P_\omega(\g^1\subset\sp(V))=\P_\omega(\g^2\subset\sp(V))=0,$$
  unless the complexification of $\g\subset\sp(V)$
 coincides with $$\so(n,\Co)\oplus\sl(2,\Co)\subset\sp(2n,\Co).$$
   For example, for $\so(n,\Real)\oplus\sl(2,\Real)\subset\sp(2n,\Real)$ it holds  $$\P_\omega(\sl(2,\Real)\subset\sp(V))=0,\quad
\P_\omega(\so(n,\Real)\subset\so(V))=(\so(n,\Co)\subset\sl(n,\Co))^{[1]}.$$

Let $\g\subset\gl(L)$ be an irreducible subalgebra. Then $\g$ is a
weakly-irreducible subalgebra of $\sp(L\oplus L^*)=\sp(2m,\Real)$,
where $m=\dim L$. Let $\omega$ be the natural symplectic form
$L\oplus L^*$. It is easy to see that $P\in\P_\omega(\g)$ if and
only if $$\pr_{\gl(L)}\circ P|_{L}\in
(\g\subset\gl(L))^{[1]},\quad \pr_{\gl(L^*)}\circ
P|_{L}\in(\g\subset\gl(L^*))^{[1]}.$$ Thus if
$\g\subset\sp(L\oplus L^*)$ is a skew-symmetric weak-Berger
subalgebra, then $(\g\subset\gl(L))^{[1]}\neq \{0\}$ and
$\g\subset\gl(L)$ is given in \cite{Odd}.

\section{Structure of the spaces $\R(\g)$}\label{secR(g)}

Consider a subalgebra $\g\subset\osp(p,q|2m)=\osp(V)$ and describe
the space $\R(\g)$. Denote the supersymmetric metric on $V$ by
$g$. It can be represented as the sum $g=\eta+\omega$, where
$\eta$ is a pseudo-Euclidean metric of signature $(p,q)$ on
$V_{\bar 0}=\Real^{p,q}$ and $\omega$ is a symplectic form on $\Pi
V_{\bar 1}=\Real^{2m}$. We identify $\osp(V)$ with $\Lambda^2V$,
the element $X\wedge Y\in\osp(V)$ is defined by
$$(X\wedge Y)Z=(-1)^{|Y||Z|}g(X,Z)Y-(-1)^{(|Y|+|Z|)|X|}g(Y,Z)X,$$
where $X,Y,Z\in V$ are homogeneous. Note that $\so(p,q)\simeq
\Lambda^2 V_{\bar 0}$ and $\sp(2m,\Co)\simeq \Lambda^2 V_{\bar
1}=\odot^2\Pi V_{\bar 1}.$
  From the Bianchi super identity it
follows that any $R\in \R(\g)$ satisfies  \begin{equation}
\label{Symsv}
g(R(X,Y)Z,W)=(-1)^{(|X|+|Y|)(|Z|+|W|)}g(R(Z,W)X,Y)\end{equation}
for all homogeneous $X,Y,Z,W\in V$. This means that
$R:\Lambda^2V\to\g\subset\Lambda^2V$ is a supersymmetric map. In
particular, $R$ is zero on the orthogonal complement
 $\g^\perp\subset\Lambda^2V$.

First consider that space $\R(\g)_{\bar 0}$. Let $R\in(\Lambda^2
V^*\otimes \g)_{\bar 0}$. Define the following maps:
\begin{align*}A&=\pr_{\so(p,q)}\g_{\bar 0}\circ R|_{\Lambda^2 V_{\bar 0}\oplus \Lambda^2 V_{\bar 1}}:\Lambda^2 V_{\bar 0}\oplus
\Lambda^2 V_{\bar 1}\to \pr_{\so(p,q)}\g_{\bar 0},\\ B&=\pr_{\Real^{2m*}\otimes\Real^{p+q}}\g_{\bar 1}\circ R|_{V_{\bar 0}\otimes V_{\bar
1}}:V_{\bar 0}\otimes V_{\bar 1}\to \pr_{\Real^{2m*}\otimes\Real^{p+q}}\g_{\bar 1},\\
C&=\pr_{\Real^{p+q\,*}\otimes\Real^{2m}}\g_{\bar 1}\circ R|_{V_{\bar 0}\otimes V_{\bar 1}}:V_{\bar 0}\otimes V_{\bar 1}\to
\pr_{\Real^{p+q\,*}\otimes\Real^{2m}}\g_{\bar 1},\\ D&=\pr_{\sp(2m,\Real)}\g_{\bar 0}\circ R|_{\Lambda^2 V_{\bar 0}\oplus \Lambda^2 V_{\bar
1}}:\Lambda^2 V_{\bar 0}\oplus \Lambda^2 V_{\bar 1}\to \pr_{\sp(2m,\Real)}\g_{\bar 0}.\end{align*} In the definition of $B$ and $C$
we used the inclusion
$$\osp(p,q|2m)\subset\gl(p+q|2m,\Real)=\gl(p+q,\Real)\oplus\gl(2m,\Real)\oplus\Real^{2m*}\otimes\Real^{p+q}\oplus
\Real^{p+q\,*}\otimes\Real^{2m}.$$ Since $R$ takes values in $\g\subset\osp(p,q|2m)$, we obtain
\begin{equation}\label{BC}\omega(C(x,\xi)y,\delta)=-\eta(y,B(x,\xi)\delta)\end{equation}
 for all $x,y\in V_{\bar 0}$ and $\xi,\delta\in V_{\bar 1}$, i.e. the maps $B$ and $C$ define each other.
  Extend the maps $A,B,C,D$ to $\Lambda^2 V$ mapping the natural complements to zero.
  Then $R=A+B+C+D$. In the matrix form we may write \begin{align*}R(x,y)&=-R(y,x)=\left(\begin{array}{cc} A(x,y)&0\\0&D(x,y)\end{array}\right),\quad
  R(\xi,\delta)=R(\delta,\xi)=\left(\begin{array}{cc} A(\xi,\delta)&0\\0&D(\xi,\delta)\end{array}\right),\\
  R(x,\xi)&=-R(\xi,x)=\left(\begin{array}{cc} 0&B(x,\xi)\\C(x,\xi)&0\end{array}\right),\end{align*}
  where $x,y\in V_{\bar 0}$ and $\xi,\delta\in V_{\bar 1}$.

  Writing down the Bianchi identity, we get that $R\in\R(\g)_{\bar 0}$ if and only if the
   following conditions hold: $A|_{\Lambda^2 V_{\bar 0}}\in\R(\pr_{\so(p,q)}\g_{\bar 0})$,
   $D|_{\Lambda^2 V_{\bar 1}}\in\bar R(\pr_{\sp(2m,\Real)}\g_{\bar
   0})$,
\begin{align}D(x,y)\xi+C(y,\xi)x+C(\xi,x)y&=0, \label{DC}\\
A(\xi,\delta)x-B(\delta,x)\xi+B(x,\xi)\delta&=0\label{AB}\end{align}
for all $x,y\in V_{\bar 0}$ and $\xi,\delta\in V_{\bar 1}$.

Suppose that $R\in\R(\g)_{\bar 0}$. Using \eqref{Symsv}, we get
\begin{align} \omega(D(x,y)\xi,\delta)&=\eta(A(\xi,\delta)x,y),\\
\omega(C(x,\xi)y,\delta)&=-\omega(C(y,\delta)x,\xi),\quad
\eta(B(x,\xi)\delta,y)=-\eta(B(y,\delta)\xi,x)
\label{SymSv2}\end{align} for all $x,y\in V_{\bar 0}$ and
$\xi,\delta\in V_{\bar 1}$.  In particular, we see that
$A|_{\Lambda^2 V_{\bar 1}}$ and $D|_{\Lambda^2 V_{\bar 0}}$ define
each other. The meanings of the restrictions \eqref{DC} and
\eqref{AB} on $A|_{\Lambda^2 V_{\bar 1}}$ and $D|_{\Lambda^2
V_{\bar 0}}$ are not so clear. On the other hand, if
representation of $\g_{\bar 0}$ is diagonal in $V_{\bar 0}\oplus
V_{\bar 1}$ (by this we mean that  the  both representations of
$\g_{\bar 0}$ on $V_{\bar 0}$ and $V_{\bar 1}$ are faithful), then
$A|_{\Lambda^2 V_{\bar 1}}$ and $D|_{\Lambda^2 V_{\bar 0}}$ are
given by $D|_{\Lambda^2 V_{\bar 1}}$ and $A|_{\Lambda^2 V_{\bar
0}}$, respectively. We will use this in many situations.

We turn now to the space  $\R(\g)_{\bar 1}$. Let $R\in(\Lambda^2
V^*\otimes \g)_{\bar 1}$. Define the following maps:
\begin{align*}A&=\pr_{\so(p,q)}\g_{\bar 0}\circ R|_{V_{\bar 0}\otimes V_{\bar 1}}:V_{\bar 0}\otimes V_{\bar 1}\to
\pr_{\so(p,q)}\g_{\bar 0},\\ B&=\pr_{\Real^{2m*}\otimes\Real^{p+q}}\g_{\bar 1}\circ R|_{\Lambda^2 V_{\bar 0}\oplus \Lambda^2 V_{\bar
1}}:\Lambda^2 V_{\bar 0}\oplus \Lambda^2 V_{\bar 1}\to \pr_{\Real^{2m*}\otimes\Real^{p+q}}\g_{\bar 1},\\
C&=\pr_{\Real^{p+q\,*}\otimes\Real^{2m}}\g_{\bar 1}\circ R|_{\Lambda^2 V_{\bar 0}\oplus \Lambda^2 V_{\bar 1}}:\Lambda^2 V_{\bar
0}\oplus \Lambda^2 V_{\bar 1}\to \pr_{\Real^{p+q\,*}\otimes\Real^{2m}}\g_{\bar 1},\\ D&=\pr_{\sp(2m,\Real)}\g_{\bar 0}\circ R|_{V_{\bar
0}\otimes V_{\bar 1}}:V_{\bar 0}\otimes V_{\bar 1}\to \pr_{\sp(2m,\Real)}\g_{\bar 0}.\end{align*}

Since $R$ takes values in $\g\subset\osp(p,q|2m)$, we obtain
\begin{equation}\label{BC1}\omega(C(x,y)z,\xi)=-\eta(z,B(x,y)\xi),\quad \omega(C(\xi,\delta)z,\theta)=-\eta(z,B(\xi,\delta)\theta)\end{equation}
for all $x,y,z\in V_{\bar 0}$ and $\xi,\delta,\theta\in V_{\bar 1}$. Thus the maps $B$ and $C$ define each other. Extend
the maps $A,B,C,D$ to $\Lambda^2V$ mapping the natural complements to zero. Then $R=A+B+C+D$. In the matrix form we may
write \begin{align*}R(x,y)&=-R(y,x)=\left(\begin{array}{cc}0&B(x,y)\\C(x,y)&0 \end{array}\right),\quad
  R(\xi,\delta)=R(\delta,\xi)=\left(\begin{array}{cc} 0&B(\xi,\delta)\\C(\xi,\delta)&0\end{array}\right),\\
  R(x,\xi)&=-R(\xi,x)=\left(\begin{array}{cc}A(x,\xi)&0\\0&D(x,\xi) \end{array}\right),\end{align*}
  where $x,y\in V_{\bar 0}$ and $\xi,\delta\in V_{\bar 1}$. Writing down the Bianchi
identity, we get that $R\in\R(\g)_{\bar 1}$ if and only if the following conditions hold:
\begin{align}B(x,y)z+B(y,z)x+B(z,x)y&=0,\label{B1}\\
C(\xi,\delta)\theta+C(\delta,\theta)\xi+C(\theta,\xi)\delta&=0,\label{C1}\\ B(x,y)\xi+A(y,\xi)x+A(\xi,x)y&=0,
\label{BA1}\\ C(\xi,\delta)x-D(\delta,x)\xi+D(x,\xi)\delta&=0\label{CD1}\end{align} for all $x,y,z\in V_{\bar 0}$ and
$\xi,\delta,\theta\in V_{\bar 1}$. Let us fix $\xi\in V_{\bar 1}$. Using \eqref{B1}, we get
$$\eta(B(x,y)\xi,z)+\eta(B(y,z)\xi,x)+\eta(B(z,x)\xi,y)=0$$ for all $x,y,z\in V_{\bar 0}$.  Using this and \eqref{BA1}, we
conclude that for each fixed $\xi\in V_{\bar 1}$ it holds
$R(\cdot,\xi)\in\P_\eta(\pr_{\so(p,q)}\g_{\bar 0})$. Similarly,
for each $x\in V_{\bar 0}$ it holds
$R(\cdot,x)\in\P_\omega(\pr_{\sp(2m,\Real)}\g_{\bar 0})$. This
will be extremely useful especially in the case when the
representation of $\g_{\bar 0}$ or of some ideal of $\g_{\bar 0}$
is diagonal in $V_{\bar 0}\oplus V_{\bar 1}$.

 In \cite{EE} it is
shown that \begin{equation}\label{RospEE}\R(\osp(p,q|2m))\simeq
\odot^2(\Lambda^2 V)/\Lambda^4 V.\end{equation} The
$\osp(p,q|2m)$-supermodule $\R(\osp(p,q|2m))$ is decomposed into
the direct sum of three irreducible components. This generalizes
the well-known decomposition of the $\so(p,q)$-module
$\R(\so(p,q))$ \cite{Al1} and defines the decomposition of the
elements $R\in\R(\osp(p,q|2m))$ into the Weyl tensor, the
trace-free part of the Ricci tensor and the scalar curvature.

Let us compare \eqref{RospEE} with the above description. For that
we consider $\R(\osp(p,q|2m))$ as an
$\so(p,q)\oplus\sp(2m,\Real)$-module. It holds \begin{align*}
\Lambda^2 V=&\Lambda^2 V_{\bar 0}\oplus (V_{\bar 0}\otimes V_{\bar
1})\oplus \Lambda^2 V_{\bar 1},\\ \odot^2( \Lambda^2
V)=&\odot^2(\Lambda^2 V_{\bar 0})\bigoplus \odot^2(V_{\bar
0}\otimes V_{\bar 1})\bigoplus \odot^2(\Lambda^2 V_{\bar
1})\bigoplus\Lambda^2 V_{\bar 0}\otimes(V_{\bar 0}\otimes V_{\bar
1})\\&\bigoplus\Lambda^2 V_{\bar 0}\otimes\Lambda^2 V_{\bar
1}\bigoplus  (V_{\bar 0}\otimes V_{\bar 1})\otimes\Lambda^2
V_{\bar 1},\\
\Lambda^4V=&\Lambda^4V_{\bar 0}\bigoplus\Lambda^3V_{\bar 0}\otimes
V_{\bar 1}\bigoplus \Lambda^2V_{\bar 0}\otimes \Lambda^2V_{\bar
1}\bigoplus V_{\bar 0}\otimes \Lambda^3V_{\bar 1} \bigoplus
\Lambda^4V_{\bar 1}.\end{align*} This implies
\begin{align*}\R(\osp(p,q|2m))_{\bar 0}\simeq& (\odot^2(\Lambda^2
V)/\Lambda^4 V)_{\bar 0}=\odot^2(\Lambda^2 V_{\bar 0})/\Lambda^4
V_{\bar 0}\\&\bigoplus\odot^2(\Lambda^2 V_{\bar 1})/\Lambda^4
V_{\bar 1}\bigoplus\odot^2(V_{\bar 0}\oplus V_{\bar
1}).\end{align*} Note that $\R(\so(p,q))\simeq\odot^2(\Lambda^2
V_{\bar 0})/\Lambda^4 V_{\bar 0}$ \cite{Al1}. Similarly,
$$\RR(\sp(2m,\Real))\simeq\odot^2(\odot^2 \Pi V_{\bar
1})/\odot^4\Pi  V_{\bar 1}=\odot^2(\Lambda^2 V_{\bar 1})/\Lambda^4
V_{\bar 1}.$$ We conclude that
$$\R(\osp(p,q|2m))_{\bar
0}\simeq\R(\so(p,q))\oplus\RR(\sp(2m,\Real))\oplus \odot^2(V_{\bar
0}\oplus V_{\bar 1}).$$ Let us describe this isomorphism. The
inclusions
$\R(\so(p,q)),\RR(\sp(2m,\Real))\subset\R(\osp(p,q|2m))_{\bar 0}$
are obvious. Let $B\in\odot^2(V_{\bar 0}\oplus V_{\bar
1})=\Lambda^2(V_{\bar 0}\oplus \Pi V_{\bar 1})$, i.e. $B$ is a
skew-symmetric endomorphism of $V_{\bar 0}\oplus \Pi V_{\bar 1}$
with respect to the skew-symmetric form $\eta\otimes \omega$, and
$B$ satisfies $$\eta\otimes\omega( B(x,\xi),y\otimes
\delta)=-\eta\otimes\omega( B(y,\delta),x\otimes \xi).$$ Note that
this corresponds to \eqref{SymSv2}.
 The
equation \eqref{BC} defines the element $C$, the equations
\eqref{DC} and \eqref{AB} defines the restrictions $D|_{\Lambda^2
V_{\bar 0}}$ and $A|_{\Lambda^2 V_{\bar 1}}$. Put in addition
$A|_{\Lambda^2 V_{\bar 0}}=0$ and $D|_{\Lambda^2 V_{\bar 1}}=0$.
We obtain  and element $R\in\R(\osp(p,q|2m))_{\bar 0}$. This
defines the inclusion $\odot^2(V_{\bar 0}\oplus V_{\bar
1})\subset\R(\osp(p,q|2m))_{\bar 0}.$

Next, $$\R(\osp(p,q|2m))_{\bar 1}\simeq V_{\bar
1}\otimes(\Lambda^2 V_{\bar 0}\otimes V_{\bar 0})/\Lambda^3
V_{\bar 0}\bigoplus V_{\bar 0}\otimes(\Lambda^2 V_{\bar 1}\otimes
V_{\bar 1})/\Lambda^3 V_{\bar 1}.$$ It holds
$\P_\eta(\so(p,q))\simeq(\Lambda^2 V_{\bar 0}\otimes V_{\bar
0})/\Lambda^3V_{\bar 0}$ \cite{onecomp}. Similarly,
$\P_\omega(\sp(2m,\Real))\simeq(\odot^2 \Pi V_{\bar 1}\otimes \Pi
V_{\bar 1})/\odot^3\Pi V_{\bar 1}.$ We obtain
$$\R(\osp(p,q|2m))_{\bar 1}\simeq \P_\eta(\so(p,q))\otimes V_{\bar
1}\bigoplus \Pi\P_\omega(\sp(2m,\Real))\otimes V_{\bar 0}.$$ Let
$P\in\P_\eta(\so(p,q))$ and $\zeta\in V_{\bar 0}$ be fixed. Define
$R\in\R(\osp(p,q|2m))_{\bar 1}$ by putting
$$A(x,\delta)=\omega(\delta,\zeta)P(x),\, D(x,\delta)=0,\,
B(x,y)\xi=\omega(\zeta,\xi)(P(x)y-P(y)x),\, B(\xi,\delta)=0.$$ In
the same way any elements $P\otimes x\in
\P_\omega(\sp(2m,\Real))\otimes V_{\bar 0}$ define an
$R\in\R(\osp(p,q|2m))_{\bar 1}$. This gives the exact form of the
obtained isomorphism.

Information about the spaces $\R(\g)$ for some Lie superalgebras
$\g\subset\gl(n|k)$ not contained in $\osp(p,q|2m)$ can be found
in \cite{EE,Poletaeva1}.

\section{Adjoint representations of simple Lie superalgebras}

\begin{prop} Let $\g$ be a simple (real or complex) Lie superalgebra
admitting an even non-degenerate $\g$-invariant bilinear supersymmetric form, i.e. such that the adjoint representation of $\g$ is
orthosymplectic. Then $\R(\g)=\R(\g)_{\bar 0}$ is one-dimensional and it is spanned by the Lie superbrackets of $\g$.
\end{prop}

{\bf Proof.} First of all, from the Jacobi super identity it follows that $[\cdot,\cdot]\in \R(\g)_{\bar 0}$ for each
simple Lie superalgebra. Note that the representation of $\g_{\bar 0}$ is diagonal in $\g_{\bar 0}\oplus \g_{\bar 1}$ (up
to the center of $\g_{\bar 0}$, which does not play a role).

First we prove that $\R(\g)_{\bar 1}=0$. Suppose that $\g_{\bar
0}$ contains at least two simple ideals $\h_1$ and $\h_2$. Let
$R\in\R(\g)_{\bar 1}$, then for each fixed $\xi\in \g_{\bar 1}$
and any $x\in\h_1$ we have $R(x,\xi)\in\h_1$. On the other hand,
for each fixed $x\in\h_1$ we have $\pr_{\g_{\bar 0}}\circ
R(x,\cdot)\in\P_\omega(\pr_{\sp(\g_{\bar 1})}\g_{\bar 0})$, but
the  $\g_{\bar 0}$-module $\g_{\bar 1}$ is a tensor product of
irreducible representations of simple ideals in $\g_{\bar 0}$ (if
$\g$ is of type I) and it is a direct sum of two such
representations (if $\g$ is of type II). This and Section
\ref{SBa} show that $\pr_{\g_{\bar 0}}\circ R(x,\cdot)$ can not
take values in a one simple ideal of $\g_{\bar 0}$ (unless
$\g_{\bar 0}$ or its complexification coincides with
$\so(n,\Co)\oplus\sl(2,\Co)$, i.e. if $\g=\osp(p,q|2)$ or
$\g=\osp(p|2,\Co)$, these cases can be considered in the same way
as $\g=\osp(1|2m,\Real)$ below and we get $\R(\g)_{\bar 1}=0$). We
have $\pr_{\g_{\bar 0}}\circ R(x,\cdot)=0$ for all $x\in\h_1$.
Similarly, $\pr_{\g_{\bar 0}}\circ R(x,\cdot)=0$ for all
$x\in\h_2$ and  $\pr_{\g_{\bar 0}}\circ R(x,\cdot)=0$ for all
$x\in\g_{\bar 0}$, i.e. $R=0$.

Suppose that the semi-simple part of $\g_{\bar 0}$ is simple, then
$\g=\osp(1|2m,\mathbb{F})$ or $\g=\osp(2|2m,\mathbb{F})$,
$\mathbb{F}=\Real$ or $\Co$ (for other simple $\g$ such that the
semi-simple part of $\g_{\bar 0}$ is simple, the adjoint
representation of $\g$ is not orthosymplectic).

Consider the case $\g=\osp(1|2m,\mathbb{F})$, the case
$\g=\osp(2|2m,\mathbb{F})$ is similar. Since the complexification
of the adjoint representation of $\g=\osp(1|2m,\mathbb{R})$ is
irreducible,  it is enough to consider the adjoint representation
of $\g=\osp(1|2m,\mathbb{C})$. Let $R\in\R(\g)_{\bar 1}$. Then for
each $x\in\g_{\bar 0}=\sp(2m,\Co)$ it holds $\pr_{\g_{\bar
0}}\circ R(x,\cdot)\in\P_\omega(\sp(2m,\Co))$, and for each
$\xi\in\g_{\bar 1}=\Co^{2m}$ it holds $\pr_{\g_{\bar 0}}\circ
R(\cdot,\xi)\in\P_\eta(\ad_{\sp(2m,\Co)})$. That is $R(\g)_{\bar
1}$ is contained in the diagonal form in the $\sp(2m,\Co)$-module
$$(\Co^{2m}\otimes  \P_\eta(\ad_{\sp(2m,\Co)}))\oplus
(\sp(2m,\Co)\otimes \P_\omega(\sp(2m,\Co))).$$ Suppose that $m\geq
2$. In \cite{onecomp} we prove that
$\P_\eta(\ad_{\sp(2m,\Co)})\simeq\sp(2m,\Co)$ and any $P\in
\P_\eta(\ad_{\sp(2m,\Co)})$ is of the form $P(\cdot)=[x,\cdot]$,
where $x\in\sp(2m,\Co)$. Note that $\P_\omega(\sp(2m,\Co))$
contains a submodule isomorphic to $\Co^{2m}$ and any element $P$
of this module is of the form $P(\cdot)=\xi \odot\cdot$ for some
$\xi\in\Co^{2m}$, here for $\xi,\delta\in \Co^{2m}$ the element
$\xi \odot\delta\in\sp(2m,\Co)$ is defined by $$(\xi
\odot\delta)\theta=\omega(\xi,\theta)\delta+\omega(\delta,\theta)\xi.$$
We conclude that $R(\g)_{\bar 1}$ is contained in the diagonal
form in the $\sp(2m,\Co)$-module $$(\Co^{2m}\otimes
\sp(2m,\Co))\oplus (\sp(2m,\Co)\otimes \Co^{2m}).$$ Moreover for
each $R\in R(\g)_{\bar 1}$ there exist linear maps
$\varphi:\Co^{2m}\to \sp(2m,\Co)$ and
$\psi:\sp(2m,\Co)\to\Co^{2m}$ such that
$R(x,\xi)=[\varphi(\xi),x]=\psi(x)\odot\xi$ for all
$\xi\in\Co^{2m}$ and $x\in \sp(2m,\Co)$. Since $R(\g)_{\bar 1}$ is
an $\sp(2m,\Co)$-module, $\varphi$ and $\psi$ are proportional as
the elements of $\Co^{2m}\otimes \sp(2m,\Co)$. To show that the
equation $[\varphi(\xi),x]=\psi(x)\odot\xi$ for all
$\xi\in\Co^{2m}$ and $x\in \sp(2m,\Co)$ has only the trivial
solution it is enough to decompose the
 $\sp(2m,\Co)$-module  $\Co^{2m}\otimes \sp(2m,\Co)$ into the direct some of irreducible components and to check this equation for a non-zero representative of each component.
We have $$\Co^{2m}\otimes \sp(2m,\Co)=V_{3\pi_1}\oplus
V_{\pi_1+\pi_2}\oplus \Co^{2m}.$$ Let
$(\xi_{\alpha})_{\alpha=-m,...,-1,1,...,m}$ be the standard basis
of $\Co^{2m}$, i.e.
$\omega(\xi_\alpha,\xi_\beta)=\delta_{\alpha,-\beta}$. Then
$\sp(2m,\Co)$ is spanned be the elements of the form
$\xi_\alpha\odot\xi_\beta$. Let
$$\varphi=c\psi=\xi_1\otimes\xi_1\odot\xi_1\in V_{3\pi_1}.$$
Substituting to the equation $\xi=\xi_{-1}$ and
$x=\xi_{-1}\odot\xi_1$, we get
$$0=[\xi_1\odot\xi_1,\xi_{-1}\odot\xi_1].$$ On the other hand,
$$[\xi_1\odot\xi_1,\xi_{-1}\odot\xi_1]=2\xi_{1}\odot\xi_1\neq 0.$$
Hence $\varphi=c\psi=\xi_1\otimes\xi_1\odot\xi_1$ is not a
solution of the equation. Similarly, taking
$$\varphi=c\psi=\sum_\alpha\xi_\alpha\otimes\xi_{-\alpha}\odot\xi_1\in
\Co^{2m},\quad \xi=\xi_{1},\quad x=\xi_{1}\odot\xi_1,$$ we get
that
$\varphi=c\psi=\sum_\alpha\xi_\alpha\otimes\xi_{-\alpha}\odot\xi_1$
is not a solution of the equation. Finally, the
$\sp(2m,\Co)$-module  $\Co^{2m}\otimes \sp(2m,\Co)$ contains the
weight space of the weight $\pi_1+\pi_2$ of dimension 2 and this
space consists of the vectors $$c_1
\xi_1\otimes\xi_1\odot\xi_2+c_2 \xi_2\otimes\xi_1\odot\xi_1,\quad
c_1,c_2\in \Co.$$ Let $\varphi=c\psi$ be equal to such vector.
Taking $\xi=\xi_{-1}$ and $x=\xi_{-2}\odot\xi_1$, we get $c_1=0$;
taking $\xi=\xi_{-2}$ and $x=\xi_{-1}\odot\xi_2$, we get $c_2=0$.
Thus we may conclude that $\R(\g)_{\bar 1}=0$. If $m=1$, then
$\Co^{2}\otimes \sp(2,\Co)=V_{3\pi_1}\oplus \Co^{2}$. It is not
hard to see that $\P_\omega(\sp(2,\Co))=\Co^2$. And the further
proof is the same.

Next we prove that  $\R(\g)_{\bar 0}$ is one-dimensional.

Let $R\in \R(\g)_{\bar 0}$ be given as above by the linear maps $A,B,C,D$.  We have seen that $A|_{\Lambda^2V_{\bar 0}}$
(or $D|_{\Lambda^2V_{\bar 1}}$)  defines uniquely $A$ and $D$. We claim that it defines the whole $R$. Indeed, suppose
that $A=0$ and $D=0$. Let $\xi,\delta\in\g_{\bar 1}$ and $x\in\g_{\bar 0}$. We have $\xi\cdot
R(x,\delta)=[\xi,R(x,\delta)]$. If $R(x,\delta)\neq 0$, then since $\g$ is simple, there exists a $\xi$ such that
$\xi\cdot R(x,\delta)\neq 0$. On the other hand, $\xi\cdot R\in\R(\g)_{\bar 1}=0$ and we get a contradiction, this proves
the claim.

If the semi-simple part of $\g_{\bar 0}$ is simple, then
$A|_{\Lambda^2V_{\bar 0}}$ being an element in $\R(\ad_{\g_{\bar
0}})$ is proportional to the Lie brackets in $\g_{\bar 0}$
\cite{Sch}. Since $A|_{\Lambda^2V_{\bar 0}}$ defines uniquely
$R\in \R(\g)_{\bar 0}$, we get that $\R(\g)_{\bar 0}$ is
one-dimensional.

Suppose that $\g$ is of type I and the semi-simple part of
$\g_{\bar 0}$ is not simple,
 then $\bar\R(\g_{\bar
0}\subset\sp(\g_{\bar 1}))$ is one-dimensional (Section
\ref{SBa}). Hence $D|_{\Lambda^2V_{\bar 1}}$ belongs to a
one-dimensional space. Since  $D|_{\Lambda^2V_{\bar 1}}$ defines
uniquely $R\in \R(\g)_{\bar 0}$, we get that $\R(\g)_{\bar 0}$ is
one-dimensional.

Finally suppose that $\g$ is of type II and the semi-simple part
of $\g_{\bar 0}$ is not simple, i.e. $\g_{\bar
0}=\h_1\oplus\h_2\oplus\z$. Then obviously
$A|_{\h_1\otimes\h_2}=0$,
$A|_{\Lambda^2\h_1}=c_1[\cdot,\cdot]_{\h_1}$, and
$A|_{\Lambda^2\h_2}=c_2[\cdot,\cdot]_{\h_2}$ are annihilated by
$\g_{\bar 0}$. Then $D|_{\Lambda^2V_{\bar 1}}$ belongs to a
two-dimensional space annihilated by $\g_{\bar 0}$. On the other
hand, $\g_{\bar 0}$ may annihilate only a one-dimensional subspace
in  $\bar\R(\g_{\bar 0}\subset\sp(\g_{\bar 1}))$. Hence there
exists a $c\in\Co$ such that for each $R$ it holds $c_1=c\,c_2$.
Thus $\R(\g)_{\bar 0}$ is one-dimensional. The proposition is
proved. $\Box$

\section{Proof of  Theorem \ref{SBTh}}

\begin{lem} Let $\g\subset\osp(p,q|2m)$ be an irreducible subalgebra of the form \eqref{vidg}. If $\R(\g)_{\bar 1}=0$, then $\R(\g)$ is a trivial $\g$-module.\end{lem}

{\it Proof.} We have $(\g_i)_{\bar 1}\cdot\R(\g)_{\bar
0}\subset\R(\g)_{\bar 1}=0$. Since each $\g_i$ is simple of
classical type, it holds $(\g_i)_{\bar 0}=[(\g_i)_{\bar
1},(\g_i)_{\bar 1}]$.  Consequently, $(\g_i)_{\bar
0}\cdot\R(\g)_{\bar 0}=0$. Suppose that $\z\neq 0$. Since
$\osp(p,q|2m)_{\bar 0}\cap\q(2m,\Real)=0$, by the Schur Lemma
$\z=\Real J$, where $J$ is an even complex structure on $V$. It is
not hard to see that $J\cdot\R(\g)=0$. $\Box$

First we consider case by case simple real Lie superalgebras $\g$ of classical type
 (they are classified in \cite{Parker80} and we list them for the convenience in
 \cite{Odd}).
For each $\g$ we find all irreducible representations
$\g\subset\osp(p,q|2m)=\osp(V)$ such that $\g$ is a non-symmetric
Berger supersubalgebra. We explain the way of the considerations
and then give several examples demonstrating this proof.

We begin with the case when $\g_{\bar 0}$ is of the form
$\h_1\oplus\h_2\oplus\z$, where $\h_1$ and $\h_2$ are simple and
$\z$ is trivial or one-dimensional.  If $\g$ is of type I, i.e.
the $\g_{\bar 0}$-module $\g_{\bar 1}$ is irreducible, then
$\g_{\bar 1}$ is of the form $W_1\otimes W_2$,
$\h_1\subset\so(W_1)$ and $\h_2\subset \sp(W_2)$ are irreducible.
If  $\g$ is of type II, i.e. the $\g_{\bar 0}$-module $\g_{\bar
1}$ is of the form $\g_{-1}\oplus\g_1$, where $\g_{-1}$ and $\g_1$
are irreducible $\g_{\bar 0}$-modules, then there are two vector
spaces $U_1$ and $U_2$ such that $\h_1\subset\gl(U_1)$,
$\h_2\subset \gl(U_2)$ are irreducible, $\g_{-1}=U_1^*\otimes
U_2$, and $\g_{1}=U_2^*\otimes U_1$.

Consider several cases:

{\it Case a.} $\h_1$ annihilates $V_{\bar 1}$. Suppose that $\g$
is of type I.
 Since the inclusion $i:\g\hookrightarrow\osp(p,q|2m)$ is a Lie superalgebra
homomorphism, the restriction $$i|_{\g_{\bar 1}}:\g_{\bar
1}=W_1\otimes W_2\to\osp(V)_{\bar 1}=V_{\bar 0}\otimes V_{\bar
1}$$ is $\g_{\bar 0}$-equivariant. In particular, it is
$\h_1$-equivariant. Since $\h_1$ annihilates $W_2$ and $V_{\bar
1}$, we conclude that $V_{\bar 0}$ is a direct sum of
$\h_1$-submodules isomorphic to $W_1$ and of an $\h_1$-trivial
submodule. Similarly, if $\g$ is of type II, then $V_{\bar 0}$ is
a direct sum of $\h_1$-submodules isomorphic to $U_1\oplus U_1^*$
and of an $\h_1$-trivial submodule.

Under the current assumption we have three cases:

{\it Case a.1.} $\h_2$ annihilates $V_{\bar 0}$. Suppose that $\g$
is of type I.
 By the above arguments, $V_{\bar 1}$ is a direct sum of $\h_2$-submodules
isomorphic to $W_2$ and of an $\h_2$-trivial submodule. From
Sections \ref{WBa} and \ref{SBa} we read that if $V_{\bar 0}$
contains more then one $\h_1$-submodule isomorphic to $W_1$ and if
$V_{\bar 1}$ contains more then one $\h_2$-submodule isomorphic to
$W_2$, then $\P_\eta(\h_1\subset\so(V_{\bar 0}))=0$ and
$\P_\omega(\h_2\subset\so(V_{\bar 1}))=0$. Consequently,
$\R(\g)_{\bar 1}=0$ and $\g$ is symmetric.  Thus either $V_{\bar
0}$ contains exactly one $\h_1$-submodule isomorphic to $W_1$, or
$V_{\bar 1}$ contains exactly one $\h_2$-submodule isomorphic to
$W_2$. Similarly, if $\g$ is of type II, then $V_{\bar 0}$
contains exactly one $\h_1$-submodule isomorphic to $U_1\oplus
U_1^*$ or $V_{\bar 1}$ contains exactly one $\h_2$-submodule
isomorphic to $U_2\oplus U_2^*$. Next we check when such
representation of $\g$ exists. For this we may pass to the
complexification  of $\g$ and of its representation. If the
resulting representation is not irreducible, we take one of its
irreducible components. Note that the type of $\g$ may change.

Let $\g$ be of type I. Suppose, for instance, that  $V_{\bar 1}$
contains exactly one $\h_2$-submodule isomorphic to $W_2$. Since
$\g_{\bar 1}\otimes W_2$ contains only one submodule annihilated
by $\h_2$ and isomorphic to $W_1$ as the $\h_1$-module and since
the representation of $\g_{\bar 1}$ on $V$ is $\g_{\bar
0}$-equivariant, $\g$ preserves the vector supersubspace
$(\g_{\bar 1}\cdot V_{\bar 1})\oplus V_{\bar 1}\subset V$  and its
even part contains only one $\h_1$-submodule isomorphic to $W_1$.
From the irreducibility of $V$ it follows that $V_{\bar 0}$
contains exactly one $\h_1$-submodule isomorphic to $W_1$ and
$V_{\bar 1}$ contains exactly one $\h_2$-submodule isomorphic to
$W_2$. Suppose that $V$ contains a non-trivial vector $v$
annihilated by  $\g_{\bar 0}$. Then the homogeneous components of
$v$ are also annihilated by  $\g_{\bar 0}$ and we may assume that
$v$ is homogeneous. Suppose that $v\in V_{\bar 0}$. Consider the
map from $\g_{\bar 1}$ to $V_{\bar 1}$ sending $\varphi\in
\g_{\bar 1}$ to $\varphi v$. This map is $\g_{\bar 0}$-equivariant
and since the $\g_{\bar 0}$-modules $\g_{\bar 1}$ and $V_{\bar 1}$
are not isomorphic, this map is zero. This shows that $\Real
v\subset V$ is an invariant subspace and it must be trivial. Thus
we get that $V_{\bar 0}=W_1$ and $V_{\bar 1}=W_2$. Such a
representation of $\g$ is either the identity one, or it does not
exist. Let $\g$ be of type II. By the similar arguments we get
that $V_{\bar 0}=U_1$ and $V_{\bar 1}=U_2$. This happens e.g. for
the complexification of the identity representation of
$\su(p_0,q_0|p_1,q_1)$.

{\it Case a.2.} $\h_2$ annihilates $V_{\bar 1}$. In this case $\h_1\oplus\h_2$ annihilates $V_{\bar 1}$. Let $U\subset
 V_{\bar 0}$ be an irreducible $\h_1\oplus\h_2$-module. Since $U$ is not $\g$-invariant, $\g_{\bar 1}\cdot U\neq 0$. On the other
hand, $\g_{\bar 1}\cdot U\subset V_1$ is annihilated by
$\h_1\oplus\h_2$, i.e. $\g_{\bar 1}\otimes U$ contains a
one-dimensional subspace annihilated by $\h_1\oplus\h_2$. This may
happens only if $U\simeq \g_{\bar 1}$ (if $\g$ is of type I), or
if $U$ is isomorphic either to $\g_{-1}$, or to $\g_{1}$ (if $\g$
is of type II). We show that such representations do not exist.

{\it Case a.3.} The representation of $\h_2$ is diagonal in
$V_{\bar 0}\oplus V_{\bar 1}$. We may decompose $V_{\bar 0}$ as
the direct sum $V_{\bar 0}=L_1\oplus L_2\oplus L_3\oplus L_4$ such
that $\h_1\oplus \h_2$ annihilates $L_4$, $\h_1$ annihilates
$L_2$, $\h_2$ annihilates $L_1$, $L_1$ is a direct sum of
$\h_1$-submodules isomorphic to $W_1$ (resp. $U_1\oplus U_1^*$),
$L_2$ is  an $\h_2$-submodule, and $L_3$ is an $\h_1\oplus
\h_2$-submodule (such that each irreducible component of $L_3$ is
faithful for both $\h_1$ and $\h_2$). If $L_3\neq 0$, then from
Section \ref{WBa} it follows that $\P_\eta(\pr_{\so(L_1\oplus
L_2)}\g_{\bar 0})=0$ and this implies $\R(\g)_{\bar 1}=0$. Thus,
$L_3=0$. Since $\g$ is a Berger algebra, one of the following
holds:

1. There exists an $R\in\R(\g)_{\bar 0}$ such that for some $x,y\in L_2$ it holds $0\neq R(x,y)\in\h_2$. Expressing $R$ in
terms of the maps $A,B,C,D$ as above, we get $A|_{\Lambda^2L_2}\neq 0$ and $D|_{\Lambda^2 L_2}\neq 0$. Hence
$A|_{\Lambda^2 V_{\bar 1}}\neq 0$ and $D|_{\Lambda^2  V_{\bar 1}}\neq 0$. This implies that $\h_2\subset \so(L_2)$ is a
Berger subalgebra and $\h_2\subset\sp( V_{\bar 1})$ is a skew-Berger subalgebra.

2. There exists an $R\in\R(\g)_{\bar 0}$ such that for some $\xi,\delta\in V_{\bar 1}$ it holds $0\neq
R(\xi,\delta)\in\h_2$. This case is similar to Case 1 and we get the same conclusion.

3. There exists an $R\in\R(\g)_{\bar 1}$ such that for some $x\in L_2$ and $\xi\in V_{\bar 1}$ it holds $0\neq
R(x,\xi)\in\h_2$. This shows that $\h_2\subset \so(L_2)$ is a Berger subalgebra and $\h_2\subset\sp( V_{\bar 1})$ is a
skew-Berger subalgebra.

Thus $L_2$ is an irreducible $\h_2$-module or $L_2=L\oplus L^*$, where $\h_2\subset \gl(L)$ is irreducible. The same holds
for $ V_{\bar 1}$. Next, as in Case a.1, we show that from the irreducibility of $\g\subset\osp(V)$  it follows that
$L_2=0$ and we get a contradiction.

{\it Case b.} $\h_1$ annihilates $V_{\bar 0}$. This case is
analogous to Case a. Note that if $\g$ is of type I, then in this
case $\h_1\subset\sp(V_{\bar 1})$ and $V_{\bar 1}$ is a direct sum
of $\h_1$-submodules isomorphic to $W_1$  and of a $\h_1$-trivial
submodule. Since $\h_1\subset\so(W_1)$, we get that
$\h_1\subset\su(W_1)$. As above we may consider Cases  b.1, b.2,
b.3.

We are left with the following case:

{\it Case c.} The representations of $\h_1$ and $\h_2$ are diagonal in $V_{\bar 0}\oplus V_{\bar 1}$.

We may decompose $V_{\bar 0}$ as the direct sum $V_{\bar 0}=L_1\oplus L_2\oplus L_3\oplus L_4$ such that $\h_1\oplus \h_2$
annihilates $L_4$, $\h_1$ annihilates $L_2$, $\h_2$ annihilates $L_1$, $L_1$ is an $\h_1$-submodule, $L_2$ is an
$\h_2$-submodule, and $L_3$ is an $\h_1\oplus \h_2$-submodule. Let $V_{\bar 0}=L'_1\oplus L'_2\oplus L'_3\oplus L'_4$ be
the similar decomposition. Since $\R(\g)_{\bar 1}\neq 0$, we get that $P_\eta(\pr_{\so(V_{\bar 0})}\g_{\bar 0})\neq 0$ and
$P_\omega(\pr_{\sp(V_{\bar 1})}\g_{\bar 0})\neq 0$. This shows that if $L_1\neq 0$, then $L_3=0$. Moreover $L_2\neq 0$, since
the representation of $\h_1$ is diagonal in $V_{\bar 0}\oplus V_{\bar 1}$. Thus either $L_3=0$, or $L_1=0$ and $L_2=0$.
Furthermore, if $L_3\neq 0$, then $\h_1\oplus \h_2\subset \so(L_3)$ is a Berger subalgebra (which is irreducible or $L_3$
is of the form $U\oplus U^*$ and  $\h_1\oplus \h_2\subset \gl(U)$ is irreducible). Similarly, if $L_1\neq 0$ and $L_2\neq
0$, then $\h_1\subset\so(L_1)$ and $\h_2\subset\so(L_2)$ are Berger subalgebras. The same statements we get for $V_{\bar
1}$ (instead of Berger subalgebras we get skew-Berger subalgebras). If $L_1\neq 0$, $L_2\neq 0$, $L'_1\neq 0$, and
$L'_2\neq 0$ then we show in the same way as in  Case a.1 that the representation is not irreducible. If  $L_1\neq 0$,
$L_2\neq 0$, and $L'_3\neq 0$, then we show as for the adjoint representations that $\R(\g)_{\bar 1}=0$. If $L_3\neq 0$,
and $L'_3\neq 0$, then either $\h_1\oplus \h_2$ does not appear as a Berger or skew-Berger algebra, or $\h_1\oplus \h_2$
appears only as a reducible Berger and a reducible skew-Berger algebra. In the last case the representation of $\g$ is not
irreducible.

\vskip0.2cm

Note that in each case we get a decomposition of $V$ into irreducible $\g_{\bar 0}$-modules. Then we ask if such
representation of $\g$ exists, in other words, we should check if the obtained representation of $\g_{\bar 0}$ can be
extended to an irreducible representation of $\g$. This can be done by passing to the complex case.  Then we may use the
theory of representations of the complex simple Lie superalgebras
\cite{Hurni82,Hurni83,Jeugt84,Kac,KacRepr,osp,F4,G3,Dict}. Any irreducible representation $\g\subset\gl(V)$ is the
highest-weight representation $V_\Lambda$ and the weight $\Lambda$ is given by its labels on the Kac-Dynkin diagram of
$\g$. There is a way to decompose the $\g_{\bar 0}$-module $V_\Lambda$ into irreducible components. One $\g_{\bar
0}$-module $V_{\tilde\Lambda}$ is obtained directly from  $\Lambda$. Then $V_{\tilde\Lambda}$ must coincide with one of
the irreducible $\g_{\bar 0}$-modules obtained by us. The weight  $\tilde\Lambda$ defines uniquely $\Lambda$ and we need
only to check that $V_\Lambda$ consists exactly of the irreducible $\g_{\bar 0}$-modules obtained by us. We will
demonstrate this technics in the examples below.

\vskip0.2cm

\begin{ex} Let $\g$ be the real form of the  complex simple Lie superalgebra $F(4)$ with
$\g_{\bar 0}=\sl(2,\Real)\oplus\so(7,\Real)$ and $\g_{\bar
1}=\Real^2\otimes \Delta$, where $\Delta\simeq\Real^8$ is the
spinor representation of $\so(7,\Real)$.

 Case a.1. We have $V_{\bar 0}=\Delta$ and $V_{\bar 1}=\Real^2$. Note that $V_{\bar 0}\otimes V_{\bar 1}=\osp(8|2,\Real)_{\bar 1}$, hence
$[V_{\bar 0}, V_{\bar 1}]=\sl(2,\Real)\oplus\so(8,\Real)$. This shows that the representation of $\g$ on $\Real^{8|2}$ does not exist.

Case a.2. We consider it after Case a.3.

Case a.3. We have $V_{\bar 0}=\Delta\oplus L$, where $\sl(2,\Real)\subset\so(L)$ is an irreducible Berger subalgebra, and
$V_{\bar 1}=\Real^2$. As in Case a.1, such representation does not exist. We may prove it also in another way. First since
$\sl(2,\Real)\subset\so(L)$ is an irreducible Berger algebra, the only possible $L$ are $\Real^3$ and $\Real^5$. Passing
to the complexification we get that as $\g_{\bar 0}=\sl(2,\Co)\oplus\so(7,\Co)$-modules, $V_{\bar 0}=\Co^8\oplus L$, where
$\Co^8$ is the spinor representation of $\so(7,\Co)$, $L$ is either $\Co^3$, or $\Co^5$, and $V_{\bar 1}=\Pi\Co^2$. Note
that neither $\g_{\bar 1}\otimes \Co^8$, nor $\g_{\bar 1}\otimes \Co^2$ contain any of the $\g_{\bar 0}$-modules $\Co^3$
and $\Co^5$. This means that $\g_{\bar 1}\cdot(\Co^8\oplus\Pi\Co^2)\subset\Co^8\oplus\Pi\Co^2$, i.e. the vector
supersubspace $\Co^8\oplus\Pi\Co^2\subset V$ is $\g$-invariant, hence $L=0$. In fact, the method of the decomposition of a
$\g_{\bar 0}$-module $V$ into irreducible components discussed above is  founded on the fact that if $U\subset V$ is an
irreducible $\g_{\bar 0}$-module, then $\g_{\bar 1}$ takes it into some irreducible components of the tensor product
$\g_{\bar 1}\otimes U$.

Let us show how the above discussed method can be applied to our
representations. The information about the irreducible
representations of $F(4)$ can be found in \cite{F4}. Any
irreducible  representation (with the highest weight $\Lambda$) of
$F(4)$ is given by the labels $(a_1,a_2,a_3,a_4)$ on the
Kac-Dynkin diagram. Define the following number
$b=\frac{1}{3}(2a_1-3a_2-4a_3-2a_4)$.  The labels must satisfy the
conditions: $b,a_2,a_3,a_4$ are non-negative integers; if $b=0$,
then $a_1=\cdots=a_4=0$; $b\neq 1$; if $b=2$, then $a_2=a_4=0$; if
$b=3$, then $a_2=2a_4+1$. The weight $\tilde \Lambda$ is given by
the labels $(b,a_2,a_3,a_4)$ on the Dynkin diagram of the Lie
algebra $\sl(2,\Co)\oplus\so(7,\Co)$. In our case $\tilde \Lambda$
must be one of $(0,0,0,1)$, $(1,0,0,0)$, $(2,0,0,0)$, $(4,0,0,0)$.
The first two cases do not satisfy the conditions on the labels.
The third case corresponds to the adjoint representation, which is
different from our ones. In the second case $V$ contains a
$\g_{\bar 0}$-module with the highest weight $(3,1,0,0)$, while
our representations do not contain such submodule.

Coming back to Case a.2 we get that in this case the representation is given by $b=1$, $a_2=0$, $a_3=0$, and $a_4=1$. But the representations of $F(4)$ with $b=1$ do not exist.

Case b. Does not appear, since the representation of $\so(7,\Real)$ in $\Delta$ is not unitary.

Case c. Section \ref{HpR} and the papers \cite{Ber57,Odd}  show
that the only representation of
 $\sl(2,\Real)\oplus\so(7,\Real)$ as a skew-Berger
algebra is in the space $\Real^2\otimes \Delta$, and Lie algebra $\sl(2,\Real)\oplus\so(7,\Real)$ does not appear as the
Berger subalgebra of $\so(p,q)$. Thus, as we have seen, in this case $\R(\g)_{\bar 1}=0$.
\end{ex}

\begin{ex} Let $\g$ be the real form of the complex simple  Lie superalgebra
$\osp(4|2,\alpha,\Co)$ with $\g_{\bar 0}=\sl(2,\Real)\oplus\sl(2,\Real)\oplus\sl(2,\Real)$ and $\g_{\bar
1}=\Real^2\otimes\Real^2\otimes\Real^2$. Suppose that $\g\subset\osp(p,q|2m)=\osp(V)$ is an irreducible Berger
supersubalgebra. We may consider several cases.

First suppose that the representation of none of the Lie algebras
$\sl(2,\Real)$ is diagonal in $V_{\bar 0}\oplus V_{\bar 1}$. Using
the fact that the representation of $\sl(2,\Real)$ in $\Real^2$ is
symplectic and the arguments of Case a.1, we get that $ V_{\bar
1}=\Real^2$ as the $\sl(2,\Real)$-module, and $ V_{\bar
0}=\Real^2\otimes\Real^2$ as the $\sl(2,\Real)\oplus
\sl(2,\Real)$-module. To analyze such representations we turn to
the complex case. Then we may use the theory of representations of
the complex simple Lie superalgebras $\osp(4|2,\alpha,\Co)$
\cite{Jeugt84}. Any representation of $\osp(4|2,\alpha,\Co)$ is
given by the numbers $(a_2,a_1,a_3)$ (the labels on the Kac-Dynkin
diagram of $\osp(4|2,\alpha,\Co)$) such that $a_2,a_3$ and the
number $b=\frac{1}{1+\alpha}(2a_1-a_2-\alpha a_3)$ are
non-negative integers. Furthermore, if $b=0$, then
$a_1=a_2=a_3=0$; if $b=1$, then $\alpha(a_3+1)=\pm (a_2+1)$.  In
\cite{Jeugt84} it is shown that $V$ contains the following
$\sl(2,\Co)\oplus \sl(2,\Co)\oplus \sl(2,\Co)$-submodules:
$(b,a_2,a_3)$,  $(b-1,a_2\pm 1,a_3\pm 1)$, $(b-1,a_2\pm 1,a_3\mp
1)$,  $(b-2,a_2\pm 2,a_3)$,  $ (b-2,a_2,a_3\pm 2)$,  $
(b-2,a_2,a_3)$,  $ (b-3,a_2\pm 1,a_3\pm 1)$, $ (b-3,a_2\pm
1,a_3\mp 1)$,  $ (b-4,a_2,a_3)$ (the representations are given by
the labels on the Dynkin diagram of $\sl(2,\Co)\oplus
\sl(2,\Co)\oplus \sl(2,\Co)$). In our case $(b,a_2,a_3)$ is one of
$(1,0,0)$,  $(0,1,0)$,  $(0,0,1)$, $(0,1,1)$, $(1,1,0)$,
$(1,0,1)$. In the first case $\alpha =1$ and we get the identity
representation of $\osp(4|2,\Real)$; the second,  the third and
the forth  cases are not possible; in the last two cases $V$
contains the representations $(0,2,1)$ and $(0,1,2)$ that give the
contradiction. Thus the only possible representation is the
identity representation of $\osp(4|2,\Real)$.

Next we suppose that a number of the representation of the Lie algebras $\sl(2,\Real)$ are diagonal in $ V_{\bar 0}\oplus
V_{\bar 1}$. By the same arguments as in Cases a.3 and b we show that if $\g$ is a Berger superalgebra, then it is
symmetric.
\end{ex}

The representations of the simple Lie superalgebras $\g$ such that
the semi-simple part of $\g_{\bar 0}$ is simple can be considered
in the same way. The situation becomes simpler, since the
representation of $\g_{\bar 0}$ is diagonal in $ V_{\bar 0}\oplus
V_{\bar 1}$ (except for the identity representations of
$\osp(1|2m,\Real)$ and $\osp(2|2m,\Real)$). We immediately
conclude that $\pr_{\so( V_{\bar 0})}\g_{\bar 0}\subset\so(V_{\bar
0})$ is a Berger subalgebra and $\pr_{\sp( V_{\bar 1})}\g_{\bar
0}\subset\sp(V_{\bar 1})$ is a skew-Berger subalgebra.

\begin{ex} Consider the Lie superalgebra $\g=\pe(n,\Real)$. Recall that $\g_{\bar 0}=\sl(n,\Real)$ and
$\g_{\bar 1}=\odot^2\Real^n\oplus\Lambda^2\Real^{n*}$. Suppose
that $\g\subset\osp(p,q|2m)=\osp(V)$ is an irreducible
non-symmetric Berger algebra. Then $\pr_{\so( V_{\bar 0})}\g_{\bar
0}\subset\so(V_{\bar 0})$ is a Berger algebra and $\pr_{\sp(
V_{\bar 1})}\g_{\bar 0}\subset\sp(V_{\bar 1})$  is a skew-Berger
algebra. Section \ref{HpR} and the results from \cite{Ber57,Odd}
show that $ V_{\bar 0}$ and $ V_{\bar 1}$ should be contained in
the following list:\\ $\sl(n,\Real)$, $\Real^n\oplus\Real^{n*}$,
$\odot^2\Real^n\oplus\odot^2\Real^{n*}$,
$\Lambda^2\Real^{n}\oplus\Lambda^2\Real^{n*}$,\\
$\Lambda^3\Real^6$ ($n=6$), and $\Lambda^4\Real^8$ ($n=8$).\\ To
study these representations we turn to the complexification.
Suppose, for example, that $ V_{\bar 0}=\Co^n\oplus\Co^{n*}$. The
$\g_{\bar 0}$-submodule $\g_{\bar 1}\cdot\Co^n$ must coincide with
an irreducible component from the list obtaining by the
complexification of the above one. We have $$\g_{\bar
1}\otimes\Co^n =V_{3\pi_1}\oplus V_{\pi_1+\pi_2}\oplus
V_{\pi_1+\pi_{n-2}}\oplus\Co^{n*}.$$ Hence, $\g_{\bar
1}\cdot\Co^n=\Co^{n*}\subset V_{\bar 1}$. The tensor product
$\g_{\bar 1}\otimes\Co^{n*}$ does not contain $\Co^{n*}$. This
means that the vector supersubspace $\Co^{n}\oplus\Co^{n*}$ (where
$\Co^{n*}\subset V_{\bar 1}$) of $V$ is $\g$-invariant and we get
a contradiction. All the other representatives from the above list
can be considered in a similar way.
\end{ex}

Thus we conclude  that if $\g$ is a simple real Lie superalgebra
and there exists an irreducible representation
$\g\subset\osp(p,q|2m)$ such that $\g$ is a non-symmetric Berger
subalgebra, then this representation is the identity one of $\g$
and $\g$ is one of the following Lie algebras with their identity
representations: $\osp(p,q|2m)$, $\osp(p|2m,\Co)$,
$\su(p_0,q_0|p_1,q_1)$ and $\hosp(p,q|m)$. These Lie superalgebras
with their identity representations are non-symmetric Berger
superalgebras, since they contain, respectively, the subalgebras
$\so(p,q)$, $\so(p,\Co)$, $\su(p_0,q_0)$ and $\hosp(p,q)$, which
are non-symmetric Berger algebras.

\vskip0.2cm

Suppose that $\g$ is a simple real Lie superalgebra  and there
exists an irreducible representation $\g\subset\osp(p,q|2m)$ such
that $\g\oplus\z$ is a non-symmetric Berger supersubalgebra, where
$\z$ is a supersubalgebra of $\osp(p,q|2m)$ commuting with $\g$.
Since $\g$ is not contained in $\q(n,\Real)$, i.e. it does not
commute with an odd complex structure, by the Schur Lemma for
representations of Lie superalgebras, $\z$ is either $\Real J$,
where $J$ is an even complex structure, or $\z=\sp(1)$, i.e. $\z$
is spanned by an even quaternionic structure $J_1,J_2,J_3$. From
Section \ref{HpR} it follows that if $\h\subset\so(p,q)$ is an
irreducible subalgebra, then if $\h\oplus\z$ is a non-symmetric
Berger algebra, then $\h$ is a non-symmetric Berger algebra.
Similarly, from \cite{Odd} it follows that if
$\h\subset\sp(2m,\Real)$ is an irreducible subalgebra  then if
$\h\oplus\z$ is a skew-Berger algebra, then $\h$ is a skew-Berger
algebra. The same note holds for representations in $U\oplus U^*$,
where $\h\subset\sl(U)$ is irreducible. This shows that the above
method can be applied also to irreducible subalgebras
$\g\oplus\z\subset\osp(p,q|2m)$, where $\g$ is simple. We obtain
the identity representations of $\u(p_0,q_0|p_1,q_1)$,
$\hosp(p,q|m)\oplus\Real J$ and $\hosp(p,q|m)\oplus\sp(1)$. Since
$\u(p_0,q_0|p_1,q_1)$ contains $\u(p_0,q_0)$, it is a
non-symmetric Berger superalgebra. Since $$\R(\sp(p,q)\oplus\Real
J)=\R(\sp(p,q)),\quad \bar\R(\so(m,\mathbb{H})\oplus\Real
J)=\bar\R(\so(m,\mathbb{H})),$$ and $J$ acts diagonally in $
V_{\bar 0}\oplus V_{\bar 1}$, we get that
$$\R(\hosp(p,q|m)\oplus\Real J)=\R(\hosp(p,q|m)),$$ i.e.
$\hosp(p,q|m)\oplus\Real J$ is not a Berger superalgebra.
Generalizing the curvature tensor of the quaternionic projective
space (with indefinite metric) \cite{Al1}, we define the curvature
tensor $R\in \R(\hosp(p,q|m)\oplus\sp(1))$ by
$$R(X,Y)=\frac{1}{2}\sum_{\alpha=1}^3g(J_\alpha
X,Y)J_\alpha-\frac{1}{4}\big(X\wedge Y+ \sum_{\alpha=1}^3J_\alpha
X\wedge J_\alpha Y\big),$$ where $X,Y\in V$. The restriction of
$R$ to $\Lambda^2\Real^{4p,4q}$ coincides with the curvature
tensor of the quaternionic projective space, its image is not
contained in $\sp(p,q)$. This shows that
$R\not\in\R(\hosp(p,q|m))$. Thus, $$\R(\hosp(p,q|m))\neq
\R(\hosp(p,q|m)\oplus\sp(1))$$ and $\hosp(p,q|m)\oplus\sp(1)$ is a
non-symmetric Berger superalgebra.

\vskip0.2cm

Let now $\g^1\subset\gl(V^1)$ and $\g^2\subset\gl(V^2)$ be two irreducible real supersubalgebras. Consider the tensor product
of these representations $\g=\g^1\oplus\g^2\subset\gl(V^1\otimes V^2)=\gl(V)$. Suppose that $\g\subset\osp(V)$.
 Note that if the real vector superspaces
$V^1$ and $V^2$ admit  complex structures commuting, respectively, with the elements of $\g^1$ and $\g^2$, then the
representation of $\g=\g^1\oplus\g^2$ in $V^1\otimes V^2$ is reducible and we consider its representation in
$V^1\otimes_{\Co} V^2$.

For the even and odd parts of $V^1\otimes V^2$ we have
$$(V^1\otimes V^2)_{\bar 0}=V^1_{\bar 0}\otimes V^2_{\bar 0}\oplus
V^1_{\bar 1}\otimes V^2_{\bar 1},\quad (V^1\otimes V^2)_{\bar
1}=V^1_{\bar 0}\otimes V^2_{\bar 1}\oplus V^1_{\bar 1}\otimes
V^2_{\bar 0}.$$ This shows that if the even and odd parts of both
$V^1$ and $V^2$ are non-trivial, then the representation of
$\g_{\bar 0}$ is diagonal in $(V^1\otimes V^2)_{\bar 0}\oplus
(V^1\otimes V^2)_{\bar 1}$. Consequently, $\pr_{\so((V^1\otimes
V^2)_{\bar 0})}\g_{\bar 0}$ is a Berger algebra and
$\pr_{\sp((V^1\otimes V^2)_{\bar 1})}\g_{\bar 0}$ is a skew-Berger
algebra. Using the  arguments as above, it is easy to see that if
$\g$ is a Berger superalgebra, then it is symmetric. The same
works for the tensor product of several representations (we may
assume that $\g^2$ is a direct some of simple Lie superalgebras
and $V^2$ is a tensor products of irreducible representations of
these Lie superalgebras).

Next, we assume that the even and odd parts of  $V^1$ are non-trivial and $V^2$ is either purely even or purely odd. If
$V^2$ is purely odd, then $V^1\otimes V^2=\Pi V^1\otimes\Pi V^2$, where $\Pi V^2$ is purely even. Thus we may assume that
$ V^2$ is purely even, i.e. $ V^2$ is a usual vector space. Since $\g=\g^1\oplus\g^2\subset\osp(V^1\otimes V^2)$, we get
that either $\g^1\subset\osp(V^1)$ and  $\g^2\subset\so(V^2)$, or $\g^1\subset\osp^{sk}(V^1)$ and  $\g^2\subset\sp(V^2)$.

Let $V^1$ and $V^2$ be a complex vector superspace and a complex
vector space, respectively, and let $V=V^1\otimes V^2$. Let $g_1$
be a supersymmetric bilinear forms on $V^1$ and $g_2$ be a
symmetric bilinear form on $V^2$ From the results of \cite{EE} it
follows that $$\R(\sl(V^1)\oplus\sl(V^2)\oplus \Co)\simeq
V^*\otimes V^*.$$ Any $\tau\in V^*\otimes V^*$ defines the
curvature tensor $R_\tau$ by
$$R_\tau(x_1\otimes x_2,u_1\otimes u_2)=A(x_1\otimes
x_2,u_1\otimes u_2)+B(x_1\otimes x_2,u_1\otimes u_2),$$ where
$A(x_1\otimes x_2,u_1\otimes u_2)\in\sl(V^1)\oplus \Co,$
$B(x_1\otimes x_2,u_1\otimes u_2)\in\sl(V^2)\oplus \Co,$ and for
$v_1\in V^1$ and $v_2\in V^2$ we have $$A(x_1\otimes x_2,u_1\otimes u_2)v_1=
(-1)^{|v_1||u_1|} \tau
(x_1,x_2,v_1,u_2)u_1-(-1)^{(|v_1|+|u_1|)|x_1|}\tau (u_1,u_2,v_1,x_2)x_1,$$ $$B(x_1\otimes x_2,u_1\otimes u_2)v_2=\tau
(x_1,x_2,u_1,v_2)u_2-(-1)^{|u_1||x_1|}\tau (u_1,u_2,x_1,v_2)x_2.$$ In particular,
$$\tr(B(x_1\otimes x_2,u_1\otimes
u_2))=(\tau(x_1,x_2,u_1,u_2)-(-1)^{|u_1||x_1|}\tau(u_1,u_2,x_1,x_2)).$$
If $q>1$, then by the same  arguments as in
\cite{MS99,Skew-Berger} it can be shown that if
$R_\tau\in\R(\osp(V^1)\oplus\so(V^2))\subset\R(\sl(V^1)\oplus\sl(V^2)\oplus\Co)$
then it is given by $\tau\in\odot^2V^*$ such that
$$\tau(x_1,x_2,u_1,u_2)=c g_1(x_1,u_1)g_2(x_2,u_2),$$ where
$c\in\Co$. Hence,
$$\R(\osp(V^1)\oplus\so(V^2))=\R(\osp(V^1)\oplus\so(V^2))_{\bar 0}$$
is one-dimensional. Thus
$\osp(V^1)\oplus\so(V^2)\subset\osp(V^1\otimes V^2)$ is a
symmetric Berger superalgebra and $\R(\g)=0$ for any proper
supersubalgebra $\g\subset
\osp(V^1)\oplus\so(V^2)\subset\osp(V^1\otimes V^2)$. Similarly, if
$g_1$ is a super-skew-symmetric bilinear forms on $V^1$ and $g_2$
is a symplectic form on $V^2$, then
$$\R(\osp^{sk}(V^1)\oplus\sp(V^2))=\R(\osp^{sk}(V^1)\oplus\sp(V^2))_{\bar 0}=\Co
R_\tau,$$
 where $$\tau(x_1,x_2,u_1,u_2)= g_1(x_1,u_1)g_2(x_2,u_2).$$ The same holds if $V^1$ and $V^2$ are real.

Thus we are left with the cases $\g^2=\sl(2,\Real)$ and
$\g^2=\sl(2,\Co)$ (in the last case $\g^1\subset\osp^{sk}(V^1)$
admits a complex structure and we consider the representation of
$\g^1\oplus\g^2$ in $V^1\otimes_\Co V^2$). Suppose that
$\g^2=\sl(2,\Real)$. We have $\g^1\subset\osp^{sk}(V^1)$. Since
there are no reductive Lie algebras $\h$ such that
$\h\oplus\sl(2,\Real)$ appears both as a Berger subalgebra of
$\so(p,q)$ and as a skew-Berger subalgebra of $\sp(2m,\Co)$, there
are no ideal in $\g^1_{\bar 0}$ that acts diagonally in $V^1_{\bar
0}\oplus V^1_{\bar 1}$. Hence $\g^1\subset\osp^{sk}(V^1)$ is the
identity representation of the Lie superalgebra $\g^1$ (this
representation must exist). Section \ref{HpR}, the results of
\cite{Odd} and the condition $\R(\g)_{\bar 1}\neq 0$ imply
$\g^1=\osp^{sk}(V^1)=\osp^{sk}(2m|r,s)$. Similarly, if
$\g^2=\sl(2,\Co)$, then $\g^1=\osp^{sk}(2m|r,\Co)$. We get the
following two algebras:
$$\osp^{sk}(2m|r,s)\oplus\sl(2,\Real)\subset\osp(\Real^{2m|r,s}\otimes
\Real^2),\,
\osp^{sk}(2m|r,\Co)\oplus\sl(2,\Co)\subset\osp(\Co^{2m|r}\otimes
\Co^2).$$ The second representation is the complexification of the
identity representation of $\hosp(r,r|m)\oplus\sp(1)$, hence the
second representation gives us a non-symmetric Berger
superalgebra. The complexification of the first representation is
the second one, hence the first representation gives us a
non-symmetric Berger superalgebra.

The theorem is proved. $\Box$

{\bf Acknowledgments.} I am grateful to D.\,V.\,Alekseevsky for
useful discussions on the topic of this paper.

\bibliographystyle{unsrt}

\end{document}